\documentclass{amsart}

\usepackage{fullpage}
\usepackage{amsfonts}
\usepackage{amssymb}
\usepackage{yhmath}

\input xy
\xyoption{matrix}\xyoption{arrow}\xyoption{curve}

\begin{document}

\renewcommand{\theequation}{\thesection.\arabic{equation}}


\newcommand{\Hom}{\mathrm{Hom}}
\newcommand{\stHom}{\underline{\mathrm{Hom}}}
\newcommand{\Ext}{\mathrm{Ext}}
\newcommand{\Tor}{\mathrm{Tor}}
\newcommand{\HH}{\mathrm{HH}}
\newcommand{\Endo}{\mathrm{End}}
\newcommand{\stEnd}{\mathrm{\underline{End}}}
\newcommand{\Tr}{\mathrm{Tr}}


\newcommand{\coker}{\mathrm{coker}}
\newcommand{\aut}{\mathrm{Aut}}
\newcommand{\op}{\mathrm{op}}
\newcommand{\add}{\mathrm{add}}
\newcommand{\ind}{\mathrm{ind}}
\newcommand{\rad}{\mathrm{rad}}
\newcommand{\soc}{\mathrm{soc}}
\newcommand{\ann}{\mathrm{ann}}
\newcommand{\im}{\mathrm{im}}
\newcommand{\chr}{\mathrm{char}}


\newcommand{\rmod}{\mbox{mod-}}
\newcommand{\Rmod}{\mbox{Mod-}}
\newcommand{\lmod}{\mbox{-mod}}
\newcommand{\lMod}{\mbox{-Mod}}
\newcommand{\stmod}{\mbox{\underline{mod}-}}
\newcommand{\stlmod}{\mbox{-\underline{mod}}}

\newcommand{\gmod}[1]{\mbox{mod}_{#1}\mbox{-}}
\newcommand{\gMod}[1]{\mbox{Mod}_{#1}\mbox{-}}
\newcommand{\Bimod}[1]{\mathrm{Bimod}_{#1}\mbox{-}}

\newcommand{\proj}{\mbox{proj-}}
\newcommand{\lproj}{\mbox{-proj}}
\newcommand{\Proj}{\mbox{Proj-}}
\newcommand{\inj}{\mbox{inj-}}
\newcommand{\coh}{\mbox{coh-}}


\newcommand{\und}[1]{\underline{#1}}
\newcommand{\gen}[1]{\langle #1 \rangle}
\newcommand{\floor}[1]{\lfloor #1 \rfloor}
\newcommand{\ceil}[1]{\lceil #1 \rceil}
\newcommand{\bnc}[2]{\scriptsize \left( \begin{array}{c} #1 \\ #2 \end{array} \right)}
\newcommand{\bimo}[1]{{}_{#1}#1_{#1}}
\newcommand{\C}{\mathcal{C}}
\newcommand{\T}{\mathcal{T}}
\newcommand{\D}{\mathcal{D}}
\newcommand{\ses}[5]{\ensuremath{0 \rightarrow #1 \stackrel{#4}{\longrightarrow} 
#2 \stackrel{#5}{\longrightarrow} #3 \rightarrow 0}}


\newtheorem{therm}{Theorem}[section]
\newtheorem{defin}[therm]{Definition}
\newtheorem{propos}[therm]{Proposition}
\newtheorem{lemma}[therm]{Lemma}
\newtheorem{coro}[therm]{Corollary}

\title{Resolutions of mesh algebras: periodicity and Calabi-Yau dimensions}
\author{Alex Dugas}
\address{Department of Mathematics, University of the Pacific, 3601 Pacific Ave, Stockton, CA 95211, USA}
\email{adugas@pacific.edu}


\begin{abstract}  A triangulated category $(\mathcal{T},\Sigma)$ is said to be Calabi-Yau of dimension $d$ if $\Sigma^d$ is a Serre functor.  We determine which stable module categories of self-injective algebras $\Lambda$ of finite type are Calabi-Yau and compute their Calabi-Yau dimensions, correcting errors in previous work.  We first show that the Calabi-Yau property of $\stmod \Lambda$ can be detected in the minimal projective resolution of the stable Auslander algebra $\Gamma$ of $\Lambda$, over its enveloping algebra.  We then describe the beginning of such a minimal resolution for any mesh algebra of a stable translation quiver and apply covering theory to relate these minimal resolutions to those of the (generalized) preprojective algebras of Dynkin graphs.  For representation-finite self-injective algebras of torsion order $t=1$, we obtain a complete description of their stable Calabi-Yau properties, but only partial results for those algebras of torsion order $t=2$.  We also obtain some new information about the periods of the representation-finite self-injective algebras of torsion order $t>1$.  Finally, we describe how these questions can also be approached by realizing the stable categories of representation-finite self-injective algebras as orbit categories of the bounded derived categories of hereditary algebras, and illustrate this technique with several explicit computations that our previous methods left unsettled.
\end{abstract}

\maketitle
 
\section{Introduction} 
\setcounter{equation}{0}

The notion of the Calabi-Yau dimension of a triangulated category was originally introduced by Kontsevich as an abstraction of a key homological characterization of Calabi-Yau varieties.  
Recently, this concept has played an important role in the study of the cluster categories introduced by Buan, Marsh, Reineke, Reiten and Todorov \cite{BMRRT} and their generalizations.  In fact, Keller and Reiten \cite{KeRe} have illustrated that the Calabi-Yau dimension is one of several invariants that characterizes the cluster categories among triangulated categories.

The goal of this article is to calculate the Calabi-Yau dimensions of the stable module categories $\stmod \Lambda$ for standard self-injective algebras $\Lambda$ of finite representation type.  Two previous papers \cite{ErdSko, BiaSko}, by Erdmann and Skowro\'{n}ski and by Bia\l kowski and Skowro\'{n}ski, respectively, already address this question.  The common strategy of these papers is to look at isomorphisms of functors on the universal cover of the category $\ind (\stmod \Lambda)$.  However, it was soon recognized that \cite{BiaSko} failed to establish the necessary functorial isomorphisms: only showing that certain syzygy functors had the same effect as the Nakayama functor on isomorphism classes of modules.  A paper of Holm and J\o rgensen \cite{HoJo1} attempted to remedy this shortcoming, claiming that such a pair of functors in a triangulated category of finite type had to be isomorphic.  Their approach employs the covering theory of triangulated categories along the lines of \cite{TOC, Amiot}.  However, while studying $\Omega$-periodicity in stable categories of finite type in \cite{Per}, we uncovered an error in some of the values of Calabi-Yau dimensions computed in \cite{ErdSko, HoJo2}, leading to a counterexample to the main result of \cite{HoJo1}.  Namely, over the representation-finite symmetric algebra of tree class $\mathbb{D}_6$ and frequency $1/3$ (see Section 6 for terminology) we found $\Omega^3 M \cong M$ for every indecomposable nonprojective module $M$, while $\Omega^3$ is not isomorphic to the identity functor (except in characteristic $2$).  Thus, the present article aims at extending the approach used in \cite{Per} to identify the Calabi-Yau dimensions of the stable categories of all self-injective algebras of finite type.  

We begin in Section 2 with the definition of Calabi-Yau dimensions of triangulated categories, and show that the Calabi-Yau property of $\stmod \Lambda$ can be expressed in terms of the homology of the stable Auslander algebra $\Gamma$ of $\Lambda$.  In this way our primary task is reduced to finding the bimodule $D\Gamma$ among the syzygies of $\Gamma$ over its enveloping algebra $\Gamma^e$.  Provided $\Lambda$ is standard, the stable Auslander algebra $\Gamma$ coincides with the mesh algebra of the stable AR-quiver of $\Lambda$, and hence we shift our focus to mesh algebras of stable translation quivers.  In Section 3, we compile some known information on these algebras and review the definitions of the {\it $m$-fold mesh algebras} from \cite{ErdSko2}, which turn out to be precisely the finite-dimensional mesh algebras.  Next, in Section 4 we give a general description of (the beginning of) the minimal projective resolution of such a mesh algebra over its enveloping algebra, generalizing resolutions computed for (generalized) preprojective algebras in \cite{ErdSna1, DPA}.  

Unfortunately, for specific algebras the precise structure of the relevant syzygies in these resolutions is still difficult to compute directly.  Hence, we approach this problem ``from the ground up'': that is, by carrying out these computations for ``minimal'' mesh algebras (for instance, the preprojective algebras) and lifting them to Galois covers.  Section 5 describes the lifting procedure developed in \cite{Per}, which relates the graded projective resolution of a group-graded algebra $A$ to the projective resolution of the Galois covering $B$ that corresponds to the given grading of $A$.  In particular, we consider the $(A,A)$-bimodule $DA$ and the Nakayama automorphism and make more explicit some of the claims in Section 6 of \cite{PAAK}.  

With all our tools in place, we then proceed in Section 6 to analyze certain graded versions of the minimal resolutions of preprojective algebras in order to compute the stable CY-dimensions of representation-finite self-injective algebras of torsion order $t=1$.  Here, we obtain the following result, which includes corrected values for the stable CY-dimensions of self-injective algebras of tree-class $\mathbb{D}_{2n}, \mathbb{E}_7$ or $\mathbb{E}_8$.  (See Section 6 for notation and terminology.)

\vspace{2mm}
\noindent
{\bf Theorem 6.1.}  {\em Suppose $\Lambda$ is a standard self-injective algebra of type $(\Delta, f, 1)$.  Then $\stmod \Lambda$ is $d$-Calabi-Yau for some $d>0$ if and only if $(h^*_{\Delta}, fm_{\Delta}) = 1$.  The minimal such $d$ is given as follows.
\begin{enumerate}
\item If $\Delta = \mathbb{A}_1, \mathbb{D}_{2n}, \mathbb{E}_7$ or $\mathbb{E}_8$ and either $2 | f$ or $\chr (k) = 2$, then $d \equiv 1- (h^*_{\Delta})^{-1} \ (\mbox{mod}\ fm_{\Delta})$ and $0 < d \leq fm_{\Delta}$.
\item In all remaining cases $d=1+2r$ where $r \equiv -(h_{\Delta})^{-1}\ (\mbox{mod}\ fm_{\Delta})$ and $0 \leq r < fm_{\Delta}$.
\end{enumerate}}

\vspace{2mm}
\noindent
In particular, this result demonstrates that the Calabi-Yau dimension of a (standard, algebraic) triangulated $k$-category of finite type is not necessarily determined by the $k$-linear structure (i.e., the AR-quiver) of the category.  Rather, we see that the characteristic of the base field $k$ plays a role as well.  This should not really come as a surprise given the important role given to $-1$ by the axioms for triangulated categories.

Section 7 is devoted to the torsion order $2$ case, requiring analysis of the graded resolutions of generalized preprojective algebras.  However, infinitely many ``minimal'' mesh algebras would need to be considered, and thus our method yields only partial results in this case (see Propositions 7.2, 7.3, 7.4).  The single family of representation-finite self-injective algebras of torsion order $t=3$ are then addressed in Section 8.  While these algebras are never stably Calabi-Yau, we calculate their periods by computing the minimal resolution of the generalized preprojective algebra associated to the Dynkin graph $\mathbb{G}_2$.

Finally, for the sake of completeness, we conclude with a description of an alternative path to many of the same results in Section 9.  This approach -- key details of which were provided by Bernhard Keller -- incorporates Asashiba's formalism \cite{Cover} to successfully modify the original strategy pursued in \cite{BiaSko, HoJo1} via the covering theory of triangulated categories.  Not only does this method make clear what was missing from the original arguments, but it also provides a simple and direct way to calculate the Calabi-Yau dimensions and periods of most (standard, algebraic) triangulated categories of finite type.  We apply it in several examples to resolve some of the ambiguities left in our original approach; nevertheless, there remain a handful of instances where we are still unable to pinpoint the Calabi-Yau dimension or period of a stable category $\stmod \Lambda$ of finite type.
  
  I would like to thank Bernhard Keller for answering my questions and kindly permitting me to include his ideas in Section 9, as well as Peter J\o rgensen and Thorsten Holm for the stimulating email discussions that greatly helped my understanding of the approach described there.  I am also indebted to the referee for suggesting the diagram in the proof of Theorem 9.5.
  
\section{Calabi-Yau dimensions of stable categories}
\setcounter{equation}{0}

Let $k$ be a field and suppose that $\mathcal{T}$ is a $k$-linear triangulated category with suspension $\Sigma$ such that $\dim_k \mathcal{T}(x,y) < \infty$ for all objects $x$ and $y$ in $\mathcal{T}$.   A triangulated autoequivalence $[S, \eta] : \mathcal{T} \rightarrow \mathcal{T}$ with $\eta : S \Sigma \stackrel{\cong}{\rightarrow} \Sigma S$ is said to be a {\it Serre functor}, or a {\it Serre duality}, if there are natural isomorphisms 
\begin{equation} t_{x,y}: \mathcal{T}(y, Sx) \stackrel{\cong}{\longrightarrow} D \mathcal{T}(x,y)
\end{equation}
 for all $x, y \in \mathcal{T}$, where $D = \Hom_k(-,k)$.  More precisely, a Serre functor should also be compatible with the suspension $\Sigma$ of $\mathcal{T}$ in the following sense \cite{CYTCcor, GCY3}: 
\begin{equation} t_{x,y} \Sigma^{-1} (\eta_x)_* = -D\Sigma t_{\Sigma x, \Sigma y}, \end{equation}
where $(\eta_x)_*$ is shorthand for the map $\T(\Sigma y, \eta_x)$.
If $\mathcal{T}$ admits an autoequivalence $S$ satisfying (2.1), then (2.2) can be used to make $S$ into a (triangulated) Serre functor in a unique way (see the appendix of \cite{GCY3}).  By definition, $\mathcal{T}$ is {\it Calabi-Yau of CY-dimension} $d$ if it has a Serre duality $S$ and there is an isomorphism of triangulated functors $[S, \eta] \cong [\Sigma, -1_{\Sigma^2}]^d \cong [\Sigma^d, (-1)^d 1_{\Sigma^{d+1}}]$ for some integer $d$.  Such a $d$ is not necessarily unique, and we thus adopt the convention that the CY-dimension of a Calabi-Yau triangulated category is the {\it smallest positive integer} $d$ such that $S \cong \Sigma^d$ as triangulated functors.  

In practice, natural isomorphisms as in (2.1) appear in various places.  However, checking the compatibility condition (2.2) is a more difficult task, perhaps requiring deeper techniques.  Thus, in the present article, we focus solely on the {\it weak Calabi-Yau dimension} of $\mathcal{T}$, which is defined to be the smallest positive integer $d$ such that $\Sigma^d \cong S$ as additive functors.  Alternatively, it is the smallest positive integer $d$ for which we have natural isomorphisms $\mathcal{T}(y, \Sigma^d x) \cong D\mathcal{T}(x,y)$ for all $x, y \in \mathcal{T}$.  Likewise, by a Serre functor we usually mean a triangulated autoequivalence satisfying (2.1), but not necessarily (2.2).  
In fact, this usage of Serre functor and of `CY-dimension' in place of `weak CY-dimension' appears quite common in the literature predating \cite{CYTC}. In particular, the weak Calabi-Yau property is often sufficient for applications (see, eg. \cite{CYTC}).

Now suppose that $A$ is a finite-dimensional self-injective $k$-algebra and consider $\mathcal{T} = \stmod A$, which is a triangulated category with the co-syzygy functor $\Omega^{-1}$ as its suspension \cite{TCRTA}.  We define the {\it stable Calabi-Yau dimension} of $A$ to be the (weak) CY-dimension of $\stmod A$ when this category is Calabi-Yau, and we define it to be infinite otherwise.  It follows from the Auslander-Reiten formula $D \stHom_A(X,Y) \cong \Ext^1_A(Y, D\Tr X)$ and the isomorphism $D\Tr \cong \Omega^2 \nu$ that $S := \Omega \nu$ is a Serre functor for $\stmod A$, where $\nu = - \otimes_A DA$ is the Nakayama functor.  Hence, we see that $\stmod A$ has Calabi-Yau dimension $d$ if and only if 
\begin{eqnarray} \Omega^{-(d+1)} & \cong & \nu \end{eqnarray}
 as triangulated functors on $\stmod A$.  As noted above, however, we will focus on checking when (2.3) is an isomorphism of (non-triangulated) functors.  Since $\Omega^{-(d+1)}$ is induced by the functor $-\otimes_A \Omega^{-(d+1)}_{A^e}(A)$, such an isomorphism would follow from an isomorphism of bimodules $\Omega^{-(d+1)}_{A^e}(A) \cong DA$.  We remark that an isomorphism of the latter variety is the basis for the definition of a {\it Calabi-Yau Frobenius algebra} in \cite{CYFA}, although it is not clear in general whether the minimal positive integer $d$ such that $\Omega^{-(d+1)}_{A^e}(A) \cong DA$ coincides with the stable CY-dimension of $A$. 

We henceforth assume that $k$ is algebraically closed and that the algebras we consider are basic.  We may thus suppose that our algebras are presented as path algebras of quivers modulo relations. Corresponding to such a coordinatization $A = kQ/I$ for a quiver $Q$, we write $e_1, \ldots, e_n$ for the primitive idempotents associated to the vertices $Q_0$ of $Q$, and we write $Q_1$ for the set of arrows of $Q$.  We denote the simple right $A$-modules (up to isomorphism) as $S_i$ for $1 \leq i \leq n$.  As we typically work with right modules, we compose paths in $kQ$ from left to right, so that we may take $Q$ to be the bound quiver of $A$.  
 
Throughout this article, we will reserve $\Lambda$ for a standard self-injective $k$-algebra of finite representation type.  We let $M_{\Lambda}$ be the direct sum of the indecomposable nonprojective $\Lambda$-modules (one representative from each isomorphism class), and let $\Gamma = \stEnd(M_{\Lambda})$ be the stable Auslander algebra of $\Lambda$.  Since $\Lambda$ is standard, $\Gamma$ is isomorphic to the mesh algebra of the stable Auslander-Reiten quiver of $\Lambda$.  We wish to relate the stable Calabi-Yau properties of $\Lambda$ and $\Gamma$, but for later applications we work in slightly greater generality.  

Assume that $\mathcal{C}$ is an exact Krull-Schmidt Frobenius $k$-category with enough projective-injectives (for instance, $\mathcal{C}$ could be $\rmod \Lambda$ or the category of Cohen-Macaulay modules over a Gorenstein hypersurface singularity).  Then the stable category $\und{\mathcal{C}}$ is triangulated, with suspension functor $\Omega^{-1}$.  If $\mathcal{C}$ has finite representation type, we can form its stable Auslander algebra $\Gamma = \stEnd_{\mathcal{C}}(M)$ where $M$ is a representation generator for $\und{\mathcal{C}}$ (without projective summands).  It is easy to see that $\rmod \Gamma$ is again a Frobenius category; in particular, if $\und{C}$ is $\Hom$-finite then $\Gamma$ is a finite-dimensional self-injective $k$-algebra.  The following proposition relates the stable Calabi-Yau property of $\und{\mathcal{C}}$ to that of $\stmod{\Gamma}$.

\begin{propos}  Assume that $\und{\mathcal{C}}$, as above, is $\Hom$-finite and has a Serre functor $S$.  Then, for a positive integer $d$, $\und{\mathcal{C}}$ is $d$-Calabi-Yau if and only if there is an isomorphism of $(\Gamma, \Gamma)$-bimodules $$\Omega^{-3d}_{\Gamma^e}(\Gamma) \cong D \Gamma.$$  In this case the stable CY-dimension of $\Gamma$ is bounded above by $3d-1$.
\end{propos}

\noindent
{\it Proof.}  Following Sections 6 and 7 of \cite{Buch}, we have an isomorphism $L := \stHom_{\mathcal{C}}(M, \Omega^{-1} M) \cong \Omega^{-3}_{\Gamma^e}(\Gamma)$ of  $(\Gamma, \Gamma)$-bimodules\footnote{Actually, for the isomorphism to hold in this generality, we need to modify Buchweitz's argument in his proof of Theorem 1.5 (cf. 7.1 in \cite{Buch}) in order to avoid the hypothesis that the Auslander algebra $\Lambda$ has Hochschild dimension $2$.  Without this hypothesis, one still obtains an epimorphism of $(\Gamma,\Gamma)$-bimodules $ \Omega^3_{\Gamma^e}(\Gamma) \oplus  (\mbox{proj.}) \rightarrow \Tor_2^{\Lambda}(\Gamma,\Gamma) \cong L$.  Since \cite{Buch} shows that these bimodules induce isomorphic functors on $\stmod \Gamma$, the kernel of this epimorphism must be a projective bimodule by Theorem 3.1 of \cite{TEG}, assuming $k$ is perfect.  Since $\Gamma$ is self-injective, the epimorphism splits and $L \cong \Omega^3_{\Gamma^e}(\Gamma)$.  We refer the reader to the proof of Theorem 3.2 in \cite{PHART} for additional details.}, which induces further bimodule isomorphisms $L^{\otimes i} \cong \stHom_{\mathcal{C}}(M, \Omega^{-i} M) \cong \Omega^{-3i}_{\Gamma^e}(\Gamma)$ for all $i \geq 1$.  We note that the left $\Gamma$-module structure on $\Omega^{-i} M$ is induced by an isomorphism $M \stackrel{\cong}{\longrightarrow} \Omega^{-i} M$, which exists since $\Omega$ is an equivalence on $\und{\mathcal{C}}$, and different choices of such an isomorphism produce isomorphic bimodules for $L$. 

Since $S$ is a Serre functor on $\und{\mathcal{C}}$, we obtain an isomorphism $$\stHom_{\mathcal{C}}(M, S M) \cong D \stHom_{\mathcal{C}}(M, M) = D\Gamma,$$
which is an isomorphism of $(\Gamma, \Gamma)$-bimodules by the naturality condition in the definition of a Serre functor.  Clearly, a functorial isomorphism $\Omega^{-d} \cong S$ on $\und{\mathcal{C}}$ induces an isomorphism $\Omega^{-d}M \cong S M$ in $\mathcal{C}$ that respects the left action of $\Gamma$, and hence we get isomorphisms $D \Gamma \cong \stHom_{\mathcal{C}}(M, \Omega^{-d}M) \cong \Omega^{-3d}_{\Gamma^e}(\Gamma)$ of bimodules over $\Gamma$.

Conversely, a bimodule isomorphism $D\Gamma \cong \Omega^{-3d}_{\Gamma^e}(\Gamma)$ induces an isomorphism $$\stHom_{\mathcal{C}}(M, S M) \cong \stHom_{\mathcal{C}}(M, \Omega^{-d}M)$$ of  $(\Gamma,\Gamma)$-bimodules.  Equivalently, we have an isomorphism of functors $\stHom_{\mathcal{C}}(-, S M) \cong \stHom_{\mathcal{C}}(-, \Omega^{-d}M)$ from $\und{\mathcal{C}}$ to $\rmod \Gamma$.  By Yoneda's lemma, it follows that $S M \cong \Omega^{-d}M$ as left $\Gamma$-modules, and this amounts to an isomorphism of functors $S \cong \Omega^{-d}$ on $\und{\mathcal{C}}$.  $\Box$ \\

The above result is key to our treatment of the standard self-injective algebras of finite representation type.  However, in characteristic $2$, there is one family of nonstandard representation-finite self-injective algebras of type $(\mathbb{D}_{3m}, 1/3, 1)$ where $m \geq 2$.  These algebras are in fact symmetric and periodic \cite{Per}, and hence stably Calabi-Yau.  Their stable Calabi-Yau dimensions are bounded above by their periods minus $1$, and we know that the period $p$ of the nonstandard algebra of type $(\mathbb{D}_{3m}, 1/3, 1)$ satisfies  $2m-1\ |\ p\ |\ 4(2m-1)$ \cite{Per}.

\section{Mesh algebras and covering theory}
\setcounter{equation}{0}

Having reduced the problem of computing CY-dimensions of representation-finite self-injective algebras to one of understanding the syzygies of the corresponding stable Auslander algebras, we now turn toward the class of mesh algebras.  In the present article, a {\it stable translation quiver} $(\Gamma, \tau)$ will mean a locally finite quiver $\Gamma = (\Gamma_0, \Gamma_1)$, possibly with loops and multiple arrows, together with a bijection $\tau : \Gamma_0 \rightarrow \Gamma_0$ such that for each pair of vertices $x, y \in \Gamma_0$, the number of arrows in $\Gamma_1$ from from $x$ to $y$ equals the number of arrows from $\tau(y)$ to $x$.  For each pair of vertices $x,y \in \Gamma_0$, we choose a bijection $\sigma : \Gamma_1(x,y) \rightarrow \Gamma_1(\tau(y),x)$, and we call $\sigma$ a {\it polarization} of $(\Gamma, \tau)$.  As $\sigma^2$ induces a bijection between the set of arrows from $x$ to $y$ and the set of arrows from $\tau(x)$ to $\tau(y)$, we extend $\tau$ to a graph automorphism of $\Gamma$ by setting $\tau(\alpha) = \sigma^2(\alpha)$ for all $\alpha \in \Gamma_1$.  We denote the initial and terminal vertices of an arrow $\alpha \in \Gamma_1$ by $i\alpha$ and $t\alpha$ respectively.

Given a stable translation quiver $(\Gamma, \tau)$, we define the {\it mesh algebra} $k(\Gamma)$ over $k$ (or over any commutative ground ring), to be the path algebra modulo relations $k\Gamma/I$, where the ideal $I$ is generated by the {\it mesh relations} $$\left\{ \sum_{t\alpha = i} \sigma(\alpha)\alpha \right\}_{i \in \Gamma_0}.$$    Note that if $\Gamma$ is infinite, the mesh algebra $k(\Gamma)$ is a ring without $1$ but with enough local idempotents.  We may also speak about the {\it mesh category} of $(\Gamma, \tau)$, which we also denote $k(\Gamma)$.  It is just the $k$-linearization of the path category of $\Gamma$, modulo the ideal $I$ generated by the mesh relations; i.e., its object set is $\Gamma_0$ and its morphism sets $k(\Gamma)(x,y)$ are quotients of the free $k$-modules generated by all paths in $\Gamma$ from $x$ to $y$ by their submodules generated by the mesh relations.  A translation quiver automorphism of $\Gamma$ naturally induces an algebra automorphism of $k(\Gamma)$, which we will denote by the same symbol.  For example, any mesh algebra has an automorphism $\tau$, corresponding to the action of $\tau$ on its quiver.

Another common definition of a translation quiver takes $\Gamma$ to be a valued quiver without multiple arrows such that the valuation of the arrow $\sigma(\alpha)$ is $(b,a)$ whenever the valuation of $\alpha$ is $(a,b)$.  Our definition may be viewed as a special case of this one by viewing $n$ parallel arrows as a single arrow with valuation $(n,n)$.  We have chosen the present definition instead in order to have $\Gamma$ coincide with the bound quiver of the mesh algebra.

\vspace{5mm}
\noindent
{\bf Examples.} (1) Given a tree $\Delta$, we recall the definition of the stable translation quiver $\mathbb{Z}\Delta$.  We first choose an orientation of $\Delta$, making it into a directed graph (the particular choice does not matter, as two different orientations lead to isomorphic translation quivers).   The vertex set of $\mathbb{Z}\Delta$ is $\Delta_0 \times \mathbb{Z}$ and there are edges $(x,n) \stackrel{(\alpha,n)}{\longrightarrow} (y,n)$ and $(y,n-1) \stackrel{(\alpha',n)}{\longrightarrow} (x,n)$ for each arrow $x \stackrel{\alpha}{\longrightarrow} y$ of $\Delta$.  We set $\tau(x,n) = (x,n-1)$ for all $(x,n) \in \Delta_0 \times \mathbb{Z}$ and define $\sigma(\alpha,n) = (\alpha',n)$ and $\sigma(\alpha',n) = (\alpha,n-1)$ for all $\alpha \in \Delta_1$ and $n \in \mathbb{Z}$. \\

\noindent
 (2) Any tree $\Delta$ can be made into a stable translation quiver by replacing each edge $\alpha$ of $\Delta$ with two arrows $\alpha_1$ and $\alpha_2$ in opposite directions.  We then set $\tau(x) = x$ for all $x \in \Delta_0$ and have $\sigma$ swap the two arrows $\alpha_1$ and $\alpha_2$ corresponding to each edge $\alpha$ of $\Delta$.  The corresponding mesh algebra is then the preprojective algebra $P(\Delta)$.  It is finite-dimensional if and only if $\Delta$ is a (simply laced) Dynkin graph.\\

\setcounter{table}{0}
\renewcommand{\thetable}{\thesection .\arabic{table}}
\begin{table}[here] $$\begin{array}{rl}
Q_{\mathbb{B}_n}: & \xymatrix{& 1 \ar[ddr] \ar@{<.>}[dd] \ar[dl]<0.5ex> & 3 \ar[l] \ar@{<.>}[dd] & \cdots \ar[l] & 2n-5 \ar[l] \ar@{<.>}[dd] \ar[ddr] & 2n-3 \ar[l] \ar@{<.>}[dd] \\ 0 \ar[ur]<0.5ex> \ar[dr]<0.5ex> \\ & 2 \ar[ul]<0.5ex> \ar[uur] & 4 \ar[l] & \cdots \ar[l] & 2n-4 \ar[l] \ar[uur] & 2n-2 \ar[l]}\\

\\
Q_{\mathbb{C}_n}: &  \xymatrix{0 \ar[dr]<0.5ex> \ar@{<.>}[dd] \\ & 2 \ar[ul]<0.5ex> \ar[r]<0.5ex> \ar[dl]<0.5ex> & 3 \ar[r]<0.5ex> \ar[l]<0.5ex> & \cdots \ar[l]<0.5ex> \ar[r]<0.5ex> & n-1 \ar[l]<0.5ex> \ar[r]<0.5ex> & n \ar[l]<0.5ex> \\ 1 \ar[ur]<0.5ex>}\\

\\
Q_{\mathbb{F}_4}: & \xymatrix{ & & 2 \ar[dl]<0.5ex> \ar@{<.>}[dd] \ar[ddr] & 4 \ar[l] \ar@{<.>}[dd] \\ 1 \ar[r]<0.5ex> & 0 \ar[l]<0.5ex> \ar[ur]<0.5ex> \ar[dr]<0.5ex> \\ & & 3 \ar[ul]<0.5ex> \ar[uur] & 5 \ar[l]}\\

\\
Q_{\mathbb{G}_2}: & \xymatrix{  & & 2 \ar[dl]<0.5ex> \ar@{<.}[dd] \\ 1 \ar[r]<0.5ex> \ar@{<.}[urr] & 0 \ar[l]<0.5ex> \ar[ur]<0.5ex> \ar[dr]<0.5ex> \\ & & 3 \ar[ul]<0.5ex> \ar@{<.}[ull]}\\

\\
Q_{\mathbb{L}_n}: & \ \ \ \  \xymatrix{0 \ar@(dl,ul) \ar[r]<0.5ex> & 1  \ar[r]<0.5ex>  \ar[l]<0.5ex> & 2  \ar[r]<0.5ex> \ar[l]<0.5ex> & \cdots \ar[l]<0.5ex> \ar[r]<0.5ex> & n-2 \ar[l]<0.5ex> \ar[r]<0.5ex> & n-1 \ar[l]<0.5ex>} 

\end{array}$$ \caption{The quivers of the generalized preprojective algebras.} 
\end{table}

\noindent
 (3) An analogue of the preprojective algebra can also be associated to each of the non-simply laced Dynkin graphs $\mathbb{B}_n, \mathbb{C}_n, \mathbb{F}_4$ and $\mathbb{G}_2$, as well as to the generalized Dynkin graphs $\mathbb{L}_n$.  The stable translation quivers are shown in Table 3.1 with the nontrivial actions of $\tau$ illustrated by dotted lines.  Whereas work of Auslander and Reiten shows that the (finite-dimensional) preprojective algebras arise as the stable Auslander algebras of the categories of maximal Cohen-Macaulay modules over simple hypersurface singularities of dimension $d=2$, Dieterich and Wiedemann have shown that many of these generalized preprojective algebras are realized by the stable Auslander algebras of simple hypersurface singularities of dimension $1$ \cite{DW} (see also \cite{Yosh}). \\

For our purposes, the most general examples of interest are the {\it $m$-fold mesh algebras} defined in \cite{ErdSko2}.  The corresponding stable translation quivers are obtained from the $\mathbb{Z}\Delta$ with $\Delta$ Dynkin, by factoring out a weakly admissible group $G$ of automorphisms (i.e., for all $g \in G -\{1\}$ and all $x \in (\mathbb{Z}\Delta)_0$, we have $x^+ \cap (gx)^+ = \emptyset$) \cite{Die}.  Such a group $G$ is always infinite cyclic and the possible generators are classified in \cite{Amiot, ADK} for $\Delta$ Dynkin.  We now list the resulting stable translation quivers for the $m$-fold mesh algebras ($m \geq 1$):
\begin{itemize}
 \item $\Delta^{(m)} = \mathbb{Z}\Delta/\gen{\tau^m}$, for $\Delta = \mathbb{A}_n, \mathbb{D}_n, \mathbb{E}_n$.
 \item $\mathbb{B}_n^{(m)} = \mathbb{ZA}_{2n-1}/\gen{\rho \tau^m}$, where $\rho = \rho_{\mathbb{A}_{2n-1}}$ is reflection in the central horizontal line of $\mathbb{ZA}_{2n-1}$.
 \item $\mathbb{C}_n^{(m)} = \mathbb{ZD}_{n+1}/\gen{\rho \tau^m}$, where $\rho=\rho_{\mathbb{D}_{n+1}}$ is induced by the order $2$ automorphism of $\mathbb{D}_{n+1}$.
 \item $\mathbb{F}_4^{(m)} = \mathbb{ZE}_6/\gen{\rho \tau^m}$, where $\rho=\rho_{\mathbb{E}_6}$ is induced by the order $2$ automorphism of $\mathbb{E}_6$.
 \item $\mathbb{G}_2^{(m)} = \mathbb{ZD}_4/\gen{\rho \tau^m}$, where $\rho$ is induced by a fixed order $3$ automorphism of $\mathbb{D}_4$.
 \item $\mathbb{L}_n^{(m)} = \mathbb{ZA}_{2n}/\gen{\rho \tau^m}$, where $\rho=\rho_{\mathbb{A}_{2n}}$ is the automorphism given by reflecting $\mathbb{ZA}_{2n}$ in the central horizontal line and then shifting half a unit to the right.  Notice $\rho^2 = \tau^{-1}$.
 \end{itemize}
It is well known that the corresponding mesh algebras are finite-dimensional and self-injective.  In fact it is shown in \cite{PAAK} that they are periodic.  It turns out that these are all of the stable translation quivers with finite-dimensional mesh algebras.  While we posit that this has been observed elsewhere, we include a brief proof below based on classical arguments as in \cite{HPR1, HPR2}.

Finally, observe that since the mesh relations always have the same form, whenever $G$ is a weakly admissible group of automorphisms of a stable translation quiver $\Gamma$, the mesh algebra $k(\Gamma)$ will be a Galois cover of the mesh algebra $k(\Gamma/G)$ with group $G$.  Furthermore, if $H \lhd G$ are weakly admissible automorphism groups of $\Gamma$, then $G/H$ is naturally a weakly admissible automorphism group of $\Gamma/H$ and $k(\Gamma/H)$ is a Galois cover of $k(\Gamma/G)$ with group $G/H$.

We conclude this section by showing that the $m$-fold mesh algebras are precisely the finite-dimensional mesh algebras of stable translation quivers.  A very similar result is proved in Theorem 4.0.4 of \cite{Amiot}, and that proof could be adapted to our setting, but we have chosen to include a slightly more direct argument.

\begin{therm}  Suppose that $(\Gamma, \tau, \sigma)$ is a (finite) stable translation quiver (possibly with loops or multiple arrows) such that the mesh algebra $k(\Gamma)$ is finite-dimensional over $k$.  Then $k(\Gamma)$ is an $m$-fold mesh algebra.
\end{therm}

\noindent
{\it Proof.} Assume $A := k(\Gamma)$ is finite-dimensional, and for all $i,j \in \Gamma_0$ let $\Gamma_1^{ij}$ be the number of arrows in $\Gamma$ from $i$ to $j$.  We first use $\Gamma$ to construct a generalized Cartan matrix $C$ (in the sense of \cite{HPR2}) and a corresponding valued graph $\Delta$.  Let $\Delta_0$ be a set of representatives of the $\tau$-orbits in $\Gamma_0$, and for each $i \in \Delta_0$, set 
$$C_{ij} = \left\{ \begin{array}{rl} - \sum_{l \in \gen{\tau}j} \Gamma_1^{il},  & i \neq j \\ 2-\sum_{l \in \gen{\tau}i} \Gamma_1^{il}, & i = j \end{array} \right.$$
Then $C = \{C_{ij}\}_{i,j \in \Delta_0}$ is a generalized Cartan matrix: i.e.,
\begin{itemize}
\item $C_{ii} \leq 2$ for all $i \in \Delta_0$;
\item $C_{ij} \leq 0$ for all $i, j \in \Delta_0$;
\item $C_{ij} \neq 0$ if and only if $C_{ji} \neq 0$, for all $i, j \in \Delta_0$.
\end{itemize}
The first two conditions are obvious, and the third follows since if $\alpha$ is an arrow from $i$ to $\tau^r(j)$, then $\sigma^{-1}\tau^{-r}(\alpha)$ is an arrow from $j$ to $\tau^{-r-1}(i)$.  (The corresponding valued graph $\Delta$ is then obtained by including an edge from $i$ to $j$ whenever $C_{ij} \neq 0$, and weighting this edge by $(|C_{ij}|, |C_{ji}|)$.)

Next, we check that setting $d_i = \dim_k (e_iA)$ for $i \in \Delta_0$ defines a subadditive function that is not additive on $C$.  To see this we will need an observation about the minimal projective resolution of a simple $A$-module: if $S_i$ is the simple $A$-module corresponding to $i \in \Gamma_0$, then the minimal projective resolution of $S_i$ begins $$e_{\tau^{-1}i}A \longrightarrow \bigoplus_{\alpha \in \Gamma_1, i\alpha = i} e_{t\alpha} A \longrightarrow e_iA \longrightarrow S_i \rightarrow 0.$$  This follows from the description of the beginning of the minimal resolution of $A$ over $A^e$ given in the next section, which holds for an arbitrary mesh algebra.  Now from this exact sequence, noting that $d_i = d_{\tau^{-1}(i)}$, we obtain the inequality $$2d_i - \sum_{i\alpha = i} d_{t\alpha} \geq 1. $$
Now consider for any $i \in \Delta_0$, 
\begin{eqnarray*} \sum_{j \in \Delta_0} d_j C_{ij} & = & d_i(2 - \sum_{l \in \gen{\tau}i} \Gamma_1^{il}) - \sum_{j \in \Delta_0 - \{i\}} d_j \sum_{l \in \gen{\tau}j} \Gamma_1^{il} \\ & = & 2d_i - \sum_{j \in \Delta_0} d_j \sum_{l \in \gen{\tau}j} \Gamma_1^{il} \\ & = & 2d_i - \sum_{j \in \Gamma_0} d_j \Gamma_1^{ij} \\ & \geq & 1.
\end{eqnarray*}
Hence, by \cite{HPR2}, $C$ is generalized Dynkin; i.e., $\Delta$ is one of the familiar graphs $\mathbb{A}_n, \mathbb{B}_n, \mathbb{C}_n, \mathbb{D}_n, \mathbb{E}_6, \mathbb{E}_7, \mathbb{E}_8, \mathbb{F}_4, \mathbb{G}_2$ or else $\mathbb{L}_n$.  We define $\Delta'$ by 
$$\Delta' = \left\{ \begin{array}{rl} \Delta, & \mbox{if\ } \Delta = \mathbb{A}_n, \mathbb{D}_n, \mathbb{E}_6, \mathbb{E}_7, \mathbb{E}_8 \\ \mathbb{A}_{2n-1}, & \mbox{if}\ \Delta = \mathbb{B}_n \\ \mathbb{A}_{2n}, & \mbox{if}\ \Delta = \mathbb{L}_n \\ \mathbb{D}_{n+1}, & \mbox{if}\ \Delta = \mathbb{C}_n \\ \mathbb{E}_6, & \mbox{if}\ \Delta = \mathbb{F}_4 \\ \mathbb{D}_4 & \mbox{if}\ \Delta = \mathbb{G}_2. \end{array} \right.$$
We now show that we can find $\Delta'$ as a sectional tree in $\Gamma$, and use it to construct a covering morphism $\mathbb{Z}\Delta' \longrightarrow \Gamma$ as in 1.6 of \cite{ADK}.  It will then follow that $\Gamma$ is a quotient of $\mathbb{Z}\Delta'$ by a weakly admissible automorphism group. 

The vertices $\Delta_0$ of $\Delta$ correspond to the $\tau$-orbits of vertices in $\Gamma$, while the edges correspond to $\tau$-orbits of arrows between $\tau$-orbits of vertices.  If $\Delta$ is simply laced, then one easily obtains a sectional tree $\Delta$ in $\Gamma$ by choosing a representative of the $\tau$-orbit corresponding to one vertex of $\Delta$ and then choosing representatives of the $\tau$-orbits of arrows incident to this vertex, etc.  We illustrate the process in the non-simply laced case with the example of $\Delta = \mathbb{F}_4$.  Here, there are $4$ $\tau$-orbits in $\Gamma_0$, say $[1], [2], [3], [4]$ corresponding to the vertices of $\mathbb{F}_4$.  
$$\xymatrix{ [1] \ar@{-}[r] & [2] \ar@{-}[r]^{(1,2)} & [3] \ar@{-}[r] & [4]}$$
Start by choosing a vertex $x \in [1] \subset \Gamma_0$.  The edge from $[1]$ to $[2]$ implies that there is an arrow from $x$ to a vertex $y \in [2]$, and then the weighted edge from $[2]$ to $[3]$ implies that there are arrows from $y$ to two different vertices $z_1, z_2 \in [3]$.  That $z_1 \neq z_2$ follows from the fact that the label is $(1,2)$, meaning that $z_1$ can only have $1$ arrow back to the orbit $[2]$.  Finally, the edge from $[3]$ to $[4]$ implies that there are arrows from both $z_1$ and $z_2$ to vertices $w_1, w_2 \in [4]$, and as above $w_1 \neq w_2$.  
$$\xymatrix{w_1 & z_1 \ar[l] & y \ar[l] \ar[d] \ar[r] & z_2 \ar[r] & w_2 \\ & & x}$$
Also notice that, by construction, this subtree of $\Gamma$ contains no paths of the form $\sigma(\alpha)\alpha$, that is, it forms a sectional $\mathbb{E}_6$-subtree.  One proceeds similarly in the other cases.  $\Box$ \\

\section{The Minimal Resolution of a Mesh Algebra}
\setcounter{equation}{0}

In this section we seek to generalize the description of the minimal projective bimodule resolutions for the preprojective algebras of (generalized) Dynkin graphs \cite{ErdSna1, PAAK, DPA}.  We let $(Q, \tau, \sigma)$ be a stable translation quiver, where $\tau$ is the translation, and $\sigma$ is the polarization.  We continue the conventions outlined in the previous section.  We shall write $A$ for the mesh algebra $k(Q, \tau, \sigma)$ of the translation quiver: In terms of quiver and relations $A = kQ/I$ where $I$ is the mesh ideal, which is generated by the relations $\{ \sum_{t(\alpha) = j} \sigma(\alpha) \alpha \ |\ j \in Q_0\}$.  Furthermore, we assume that $A$ is finite-dimensional and self-injective (actually, self-injectivity follows from finite-dimensionality by Theorem 3.1 for instance).  In particular, we can fix a bilinear form $(-,-)$ on $A$ corresponding to an isomorphism $A \cong DA$.  

We will show that the beginning of the minimal projective resolution of $A$ over its enveloping algebra $A^e = A^{\op} \otimes_k A$ has the form
\begin{eqnarray} \ses{L \longrightarrow P_2}{P_1}{P_0 \longrightarrow A}{R}{\delta} \end{eqnarray}
where $L$ is a twisted bimodule ${}_1 A_{\mu}$ with $\mu \in \aut(A)$, $P_0 = \oplus_{i \in Q_0} (e_i \otimes e_i)A^e$, $P_1 = \oplus_{\alpha \in Q_1} (e_{i\alpha} \otimes e_{t\alpha})A^e$ and $P_2 = \oplus_{i \in Q_0} (e_{\tau i} \otimes e_i)A^e$; and where 

\begin{eqnarray} \delta(e_{i\alpha} \otimes e_{t\alpha})  =  x_{\alpha} & := & \alpha \otimes e_{t\alpha} - e_{i\alpha} \otimes \alpha, \ \ \ \mbox{for}\ \alpha \in Q_1 \\ R(e_{\tau i} \otimes e_i)  =  \sigma_i & := & \sum_{t \alpha = i} (\sigma(\alpha) \otimes e_{t\alpha} + e_{\tau i} \otimes \alpha), \ \ \mbox{for}\ i \in Q_0. 
\end{eqnarray}

Minimality and exactness of the first two terms of the resolution follow from \cite{HCFDA}.  Next, we check that $\delta R = 0$:
\begin{eqnarray*} \delta R( e_{\tau i} \otimes e_i) & = & \sum_{t \alpha = i} (\sigma(\alpha) \delta(e_{i\alpha} \otimes e_{t\alpha}) + \delta(e_{i \sigma \alpha} \otimes e_{t\sigma \alpha}) \alpha) \\ & = & \sum_{t\alpha = i} (\sigma(\alpha)\alpha \otimes e_{t\alpha} - \sigma(\alpha) \otimes \alpha + \sigma(\alpha) \otimes \alpha - e_{i\sigma \alpha} \otimes \sigma(\alpha) \alpha) \\ & = & 0.
\end{eqnarray*}

We now show that $\im (R) = \ker (\delta)$.  Observe that none of the $\sigma_i$ belong to the radical of $\ker (\delta)$.  We also know from the generators of the mesh ideal $I$ that we have a projective cover $\pi : P_2 \rightarrow \ker (\delta)$.  Thus $R = \pi h$ for some endomorphism $h$ of $P_2$, and the images of the generators of $P_2$ under $h$ cannot lie in the radical of $P_2$.  Furthermore, since these  generators are normed by different primitive idempotents $e_i$ (on the right, say), their images under $h$ will be linearly independent modulo the radical of $P_2$.  In other words, $h$ is an isomorphism, and we can identify $\pi$ with $R$.

Using the fact that $A$ must be an $m$-fold mesh algebra, it follows from \cite{DPA, DMA} that $L = \Omega^3_{A^e}(A)$ is a twisted bimodule of the form ${}_1 A_{\mu}$ for some $\mu \in \aut_k(A)$.  The rest of this section is devoted to a more detailed description of $L$, parallel to work in \cite{ErdSna1}, which can be applied to calculate the automorphism $\mu$ in particular cases.

We begin by describing a set of generators of $L$.  We first fix a basis $\mathcal{B}$ for $A$, consisting of homogeneous elements (with respect to the path-length grading of $A$) such that $\mathcal{B} = \cup_{i, j \in Q_0} e_i \mathcal{B} e_j$.  Using the bilinear form $(-,-)$ on $A$, we can define a dual basis $\mathcal{B}^*$ so that $(b_i, b_j^*) = \delta_{i,j}$ for all $b_i, b_j \in \mathcal{B}$.  Notice that $b \in e_i \mathcal{B} e_j$ if and only if $b^* \in e_j \mathcal{B}^*e_{\pi(i)}$, where $\pi$ is the Nakayama permutation of $A$ (written as a permutation of the indices $i$).  Letting $|x|$ denote the degree of $x \in \mathcal{B}$ with respect to the path length grading of $A$, we can now define
\begin{equation} \xi_i  :=  \sum_{x \in e_i \mathcal{B}} (-1)^{|x|} \tau(x) \otimes x^* \in P_2 \label{eq:xidef} \end{equation}
for each $i \in Q_0$.

\begin{lemma} For each $i \in Q_0$, we have $R(\xi_i) = 0$.
\end{lemma}

\noindent
{\it Proof.} We have \begin{eqnarray*}
R(\xi_i) & = & \sum_{x \in e_i \mathcal{B}} (-1)^{|x|} \tau(x) R(e_{\tau(tx)} \otimes e_{tx}) x^* \\ 
& = & \sum_{x \in e_i \mathcal{B}} (-1)^{|x|} \sum_{t\alpha = tx} \tau(x) (\sigma(\alpha) \otimes e_{tx} + e_{\tau(tx)} \otimes \alpha) x^* \\
& = & \sum_{\alpha \in Q_1} \sum_{d \geq 1} (-1)^{d-1} \left[ \sum_{x \in e_i \mathcal{B} e_{t\alpha}, |x| = d-1} \tau(x) \sigma(\alpha) \otimes x^* - \sum_{x \in e_i \mathcal{B} e_{\tau^{-1}(i\alpha)}, |x| = d} \tau(x) \otimes \sigma^{-1}(\alpha) x^*\right].
\end{eqnarray*}
Now we wish to show that each of the summands in the last sum vanishes.  Thus, we fix $e_i, d$ and $\alpha$, and imitating \cite{ErdSna1}, we write $$e_i \mathcal{B}_{d-1}e_{t\alpha} = \{y_1, \ldots, y_r\}\ \ \ \mathrm{and}\ \ \ e_i \mathcal{B}_d e_{\tau^{-1}(i\alpha)} = \{x_1, \ldots, x_s\}.$$  The summand enclosed in brackets above thus becomes 
\begin{eqnarray}\sum_{i = 1}^r \tau(y_i) \sigma(\alpha) \otimes y_i^* - \sum_{j=1}^s \tau(x_j) \otimes \sigma^{-1}(\alpha)x_j^*.
\end{eqnarray}
Observe that $\tau(y_i) \sigma(\alpha) = \tau(y_i \sigma^{-1}(\alpha))$ and
\begin{eqnarray}
y_i \sigma^{-1}(\alpha) & = & \sum_{j=1}^s (y_i \sigma^{-1}(\alpha), x_j^*) x_j \\
& = & \sum_{j=1}^s (y_i, \sigma^{-1}(\alpha) x_j^*)x_j.
\end{eqnarray} 
Hence $\sum_{i = 1}^r \tau(y_i) \sigma(\alpha) \otimes y_i^* = \sum_{i=1}^r \sum_{j=1}^s (y_i, \sigma^{-1}(\alpha)x_j^*) \tau(x_j) \otimes y_i^*$, and similarly $\sum_{j=1}^s \tau(x_j) \otimes \sigma^{-1}(\alpha)x_j^* = \sum_{j=1}^s \sum_{i=1}^r \tau(x_j) \otimes (y_i, \sigma^{-1}(\alpha)x_j^*)y_i^*$.  Thus, we obtain the aforementioned cancellation.  $\Box$ \\

\begin{lemma}  For each $i \in Q_0$ we have isomorphisms $$\xi_i A \cong e_{\pi i} A, \ \ \ \mathrm{and}\ \ \ \ A\xi_i \cong A e_{\tau i}.$$
Thus,  letting $\xi = \sum_{i \in Q_0} \xi_i$, we get $\xi A = \oplus_{i \in Q_0} \xi_i A \cong A_A$ and $A \xi = \oplus_{i \in Q_0} A \xi_i \cong {}_A A$.
\end{lemma}

\noindent
{\it Proof.}  By construction we have $\xi_i = e_{\tau i} \xi_i e_{\pi i}$.  Thus the natural map $e_{\pi i}A \rightarrow \xi_i A$ given by $a \mapsto \xi_i a$ is onto.  To see that it is injective, it suffices to show that $e_{\pi i}^*$ is not sent to $0$, since this element generates the socle of $e_{\pi i}A$.  The image of $e_{\pi i}^*$ is $$\xi_i e_{\pi i}^* = (-1)^{|e_i^*|} \tau(e_i^*) \otimes e_{\pi i}^* \neq 0$$ since $e_{\pi i}^*$ is annihilated by the radical of $A$ on either side.  The proof of the isomorphism $A \xi_i \cong A e_{\tau i}$ is similar. $\Box$ \\

\begin{coro}  $L = \ker(R)$ is generated by $\xi$, and is isomorphic to a twisted bimodule ${}_1 A_{\mu}$ for some $k$-algebra automorphism $\mu$ of $A$.  Moreover $\mu(e_i) = e_{\pi \tau^{-1}(i)}$ for each $i \in Q_0$.
\end{coro}

\noindent
{\it Proof.}  Since we already know $L$ is a twisted bimodule, the above lemma shows that $\xi$ is a generator.  To describe $\mu$, we use the left module isomorphism ${}_A L \rightarrow {}_A A$ induced by $\xi \mapsto 1$ to transport the bimodule structure of $L$ to $A$, obtaining ${}_1 A_{\mu}$ by definition of $\mu$.  This yields the identity $a \xi = \xi \mu(a)$ for any $a \in A$.  Now letting $a = e_i$, the equality $\xi_i = e_{\tau i} \xi_i e_{\pi i}$ easily leads to $\mu(e_i) = e_{\pi\tau^{-1}(i)}$.  $\Box$ \\

\section{Lifting the Nakayama functor}
\setcounter{equation}{0}

As shown above, deciding whether a representation-finite self-injective algebra $\Lambda$ is stably Calabi-Yau amounts to determining if some co-syzygy of its stable Auslander algebra $\Gamma$ is isomorphic to the bimodule $D\Gamma$.  Since $\Gamma$ is a mesh algebra of Dynkin type, we can use the theory of Galois covers and smash products to relate it to a preprojective algebra of a (simply laced) Dynkin graph.  In this section, we review the necessary results on covering theory, and apply them to show that whenever an algebra $B$ is a finite Galois cover of $A$, the bimodule lifting functor of \cite{Per} takes $DA$ to $DB$.  We also discuss how to lift the Nakayama automorphism  from $A$ to $B$ in this case.

In order to analyze Galois coverings of an algebra $A$ with finite Galois group $G$, we will study the corresponding $G$-grading of $A$ and the smash product algebra $A \# k[G]^*$ as in \cite{CiMa, CoMo}.  In this setting, the grading of $A$ arises from a function $\pi : Q_1 \rightarrow G$ as in \cite{Green}, whose image generates $G$.  In particular, $\deg(e_i) = e$ (the identity of $G$) for all $i$, and the Jacobson radical $J_A$ is a homogeneous ideal containing $\oplus_{g \neq e} A_g$.  Furthermore we assume that $A$ is indecomposable and the grading is such that $B$ is indecomposable as well.  We write $\gmod{G} A$ for the category of finite-dimensional $G$-graded right $A$-modules and degree-preserving morphisms.  We also write $\Bimod{G} A$ (respectively, $\Bimod{} A$) for the category of $G$-graded (respectively, all) $(A,A)$-bimodules.  For a graded $A$-module $M$ and $d \in G$, we define $M[d]$ to be the graded $A$-module given by $M[d]_g = M_{d^{-1}g}$ for all $g \in G$.

The smash product $B = A \# k[G]^*$ is the free $A$-module $\oplus_{g \in G} A p_g$ with multiplication $ap_g \cdot b p_h = ab_{gh^{-1}} p_h$, where $b_g$ denotes the degree-$g$ component of an element $b \in A$.  It is easy to see that the elements $\{e_i p_g\ |\ 1 \leq i \leq n,\ g \in G\}$ form a complete set of pairwise orthogonal primitive idempotents for $B$.  It is also well known that there is an isomorphism of categories $\gmod{G} A \cong \rmod B$.

In \cite{Per}, we defined a functor $F = F_e : \Bimod{G} A \rightarrow \Bimod{} B$ by setting $F(M) = M \otimes_A B$ as an $(A,B)$-bimodule, and extending the left action of $A$ to $B$ by the rule $$ap_g \cdot m_k p_h = \left\{ \begin{array}{cc} am_k p_h, & \mathrm{if}\ g = kh \\ 0, & \mathrm{otherwise}  \end{array} \right.$$ for $a \in A, g \in G, m = \sum_{k \in G} m_k \in M$ with $m_k \in M_k$ and $k, h \in G$.  
It is shown in \cite{Per} that $F$ is exact and takes projectives to projectives.  Furthermore, $F(M[d]) \cong {}_1 F(M)_d$, whenever $d \in Z(G)$.  We wish to study the effect of this functor on the bimodule $DA$.  Notice that $DA = \Hom_k(A,k)$ has a natural $G$-grading given by $(DA)_g = \Hom_k(A_{g^{-1}},k)$ for each $g \in G$.  Henceforth, we assume $DA$ is given this grading.

\begin{propos} There is a $(B,B)$-bimodule isomorphism $F(DA) \cong DB$.
\end{propos}

\noindent
{\it Proof.}  We define $\varphi : DA \otimes_A B \rightarrow DB$ by $$\varphi(f p_g)(a p_h) = f(a_{gh^{-1}})$$ for all $f \in DA, a \in A$ and $g, h \in G$.  We first check that $\varphi$ is a right $B$-module homomorphism.  For $a, b \in A, f \in DA$ and $g, h, l \in G$, we have 
\begin{eqnarray*} \varphi(fp_g \cdot b p_l)(a p_h) & = & \varphi(f b_{gl^{-1}} p_l)(a p_h) \\ 
& = & f b_{gl^{-1}}(a_{lh^{-1}}) \\ & = & f (b_{gl^{-1}}a_{lh^{-1}})
\end{eqnarray*}
and
 \begin{eqnarray*} (\varphi(fp_g) \cdot b p_l)(a p_h) & = & \varphi(fp_g)( b a_{lh^{-1}} p_h) \\ 
& = & f (b_{gl^{-1}} a_{lh^{-1}}).
\end{eqnarray*}
To see that $\varphi$ is a right $B$-module homomorphism, let $f_k \in (DA)_{k}$ for some $k \in G$ and consider
\begin{eqnarray*} \varphi(b p_g \cdot f_k p_l)(a p_h) & = & \left\{ \begin{array}{cc} \varphi( b f_k p_l)(a p_h), & \mathrm{if}\ g = kl \\ 0, & \mathrm{if}\ g \neq kl \end{array} \right.  \\ 
& = & \left\{ \begin{array}{cc}  f_k(a_{lh^{-1}}b), & \mathrm{if}\ g = kl \\ 0, & \mathrm{if}\ g \neq kl \end{array} \right.  \\ 
& = & \left\{ \begin{array}{cc} f_k(a b_{hg^{-1}}), & \mathrm{if}\ g = kl \\ 0, & \mathrm{if}\ g \neq kl \end{array} \right.  \\ 
& = & \varphi(f_k p_l)( ap_h b p_g) \\ & = & (bp_g \cdot \varphi(f_k p_l))(ap_h).
\end{eqnarray*}
To show that $\varphi$ is an isomorphism clearly it suffices to show that it is injective, since both modules have the same $k$-dimension.  Thus suppose that $\varphi(\sum_{h \in G}f^h p_h) = 0$ for certain $f^h \in DA$, which are not all $0$.  There exists a $g \in G$ and a homogeneous element $a = a_l \in A$ such that $f^g(a_l) \neq 0$, and then $\varphi(\sum_{h \in G} f^h p_h)(a p_{l^{-1}g}) = \sum_{h \in G} f^h(a_{hg^{-1}l}) = f^g(a_l) \neq 0$, a contradiction.  $\Box$ \\

If $A$ is self-injective, the bimodule $DA$ is isomorphic to a twisted bimodule ${}_{\nu} A_1$ where $\nu \in \aut(A)$ is called a {\it Nakayama automorphism} of $A$.  In general, $\nu$ represents a unique element in the group of outer $k$-algebra automorphisms of $A$.  We write $\pi$ for the {\it Nakayama permutation} of $A$, which is the permutation of the isomorphism classes of primitive idempotents of $A$ that is induced by $\nu$.  As a consequence of the above proposition, we recover the well-known fact that finite Galois covers of self-injective algebras are again self-injective: for, as right $B$-modules $DB_B \cong DA \otimes_A B_B \cong A \otimes_A B_B \cong B_B$.  Furthermore, if $\nu$ is a graded automorphism of $A$, we can describe the Nakayama automorphism of $B$.  In this case, the twisted bimodule ${}_{\nu} A_1$ is graded with the same grading as $A$, but may not be isomorphic to $DA$ as a graded bimodule.

\begin{lemma}  Suppose a Nakayama automorphism $\nu$ of $A$ is graded, i.e., $\nu(A_g) = A_g$ for all $g \in G$.  Then $B$ has Nakayama automorphism $\tilde{\nu}$ defined by $$\tilde{\nu}(ap_g) = \nu(a)p_{gx}$$ for some $x \in G$.  Furthermore, if $G$ is abelian, then $x$ is the unique element of $G$ such that $DA[x] \cong {}_{\nu} A_1$ as graded bimodules.
\end{lemma}

\noindent
{\it Proof.}  Since $DA \cong {}_{\nu} A_1$ as ungraded bimodules and $F(DA) \cong DB$ is indecomposable, Lemma 3.5 in \cite{Per} shows that $DB \cong {}_x F({}_{\nu} A_1)$ for some $x \in G$.   As a right module, ${}_x F({}_{\nu} A_1) \cong A \otimes_A B \cong B_B$ via the identity map $b p_h \mapsto bp_h$, and hence we can obtain $\tilde{\nu}$ by examining the left-action of $B$ on this bimodule.  For $a, b \in A$ and $g, h, k \in G$, we have $$ap_g \cdot (b_k p_h) = a p_{gx} \cdot b_k p_h = \left\{ \begin{array}{cc} \nu(a) b_k p_h & \mathrm{if}\ gx = kh \\ 0 & \mathrm{otherwise.} \end{array} \right.$$ One easily checks that this agrees with the left action of $B$ on ${}_{\tilde{\nu}} B_1$ with $\tilde{\nu}$ as in the statement of the lemma.  If $G$ is abelian, then there exists a unique $x$ such that $DA[x] \cong {}_{\nu} A_1$ as graded bimodules, since in this case graded bimodules can be identified with graded modules over the enveloping algebra $A^e$.  Thus $DB \cong F(DA) \cong F({}_{\nu}A_1[x^{-1}])\cong {}_xF({}_{\nu}A_1)$, showing that this $x$ agrees with the one used above.  $\Box$ \\

The algebra $A$ is said to be {\it symmetric} if its Nakayama automorphism $\nu$ is the identity.  In this case, the lemma shows that the Nakayama automorphism $\tilde{\nu}$ for $B$ is the automorphism $ap_g \mapsto ap_{gx}$ induced by some $x \in G$.  Although the converse does not hold in general, the following observations concerning Nakayama permutations shall suffice for our purposes.  Notice that $\tilde{\nu}(e_i p_g) = \nu(e_i) p_{gx} \cong \pi(e_i) p_{gx}$, and thus the Nakayama permutation $\tilde{\pi}$ for $B$ is given by $\tilde{\pi}(e_i p_g) = \pi(e_i) p_{gx}$.  In particular, we see that the Nakayama permutation for $B$ is induced by some automorphism $x \in G$ if and only if $\pi$ is the identity, i.e., if and only if $A$ is weakly symmetric.

\section{Preprojective algebras and the torsion order 1 case}
\setcounter{equation}{0}

We now apply our results on mesh algebras to study the standard self-injective algebras of finite representation type.  We recall that if $\Lambda$ is a standard self-injective algebra of finite type with stable AR-quiver $(\Gamma, \tau)$, the full subcategory $\ind (\stmod \Lambda)$ of indecomposable objects in the stable category $\stmod \Lambda$ is equivalent to the mesh category $k(\Gamma)$, and the stable Auslander algebra of $\Lambda$ is isomorphic to the mesh algebra $k(\Gamma)$.  Furthermore, Asashiba has shown that all stably equivalent standard self-injective algebras of finite representation type are in fact derived equivalent \cite{DECSA}.  It then follows from \cite{DCSE} that the stable AR-quiver of $\Lambda$ determines the stable category $\stmod \Lambda$ as a triangulated $k$-category and hence also determines the stable CY-dimension of $\Lambda$.

  As for the possible translation quivers that can arise as the AR-quiver of $\Lambda$, Riedtmann's classification tells us that they have the form $\mathbb{Z}{\Delta}/\Pi$ for a simply laced Dynkin graph $\Delta$ and admissible group of automorphisms $\Pi$ of the translation quiver $\mathbb{Z}\Delta$ \cite{ADK}.  Moreover, $\Pi$ is infinite-cyclic, generated by $\zeta \tau^{-r}$, where $\tau$ is the translation and $\zeta$ is an automorphism of finite order $t$.  Setting $f = r/m_{\Delta}$, the triple $(\Delta, f, t)$ is then known as the {\it type} of $\Lambda$, and $\Delta, f, t$ are the {\it tree-class, frequency} and {\it torsion order} of $\Lambda$, respectively.  Here, $m_{\Delta}$ equals $n, (2n-3), 11, 17$ or $29$ when $\Delta$ is $\mathbb{A}_n, \mathbb{D}_n, \mathbb{E}_6, \mathbb{E}_7$ or $\mathbb{E}_8$ respectively, and the {\it Coxeter number} of $\Delta$ is $h_{\Delta} = m_{\Delta} + 1$.  The possible types were classified by Riedtmann, and we refer the reader to the appendix of \cite{Asa2} for a list, along with presentations of representative algebras of each type in terms of quivers and relations.  Notice that since the type of $\Lambda$ completely specifies the AR-quiver of $\Lambda$, we will be able to describe the stable CY-dimension of $\Lambda$ solely in terms of these data.  To simplify the statement of our main result, we define
 \begin{eqnarray} h^*_{\Delta} & = & \left\{ \begin{array}{cc} h_{\Delta}, &  \Delta = A_n, D_{2n+1}, E_6, \ (n \geq 2) \\ h_{\Delta}/2, & \Delta = A_1, D_{2n}, E_7, E_8; \ (n \geq 2) \end{array} \right. \end{eqnarray}
Note that this differs from the definition in \cite{BiaSko} when $\Delta = \mathbb{A}_{4l+1}$ with $l \geq 1$.

\begin{therm} Suppose $\Lambda$ is a standard self-injective algebra of type $(\Delta, f, 1)$.  Then $\stmod \Lambda$ is $d$-Calabi-Yau for some $d>0$ if and only if $(h^*_{\Delta}, fm_{\Delta}) = 1$.  The minimal such $d$ is given as follows.
\begin{enumerate}
\item If $\Delta = \mathbb{A}_1, \mathbb{D}_{2n}, \mathbb{E}_7$ or $\mathbb{E}_8$ and either $2 | f$ or $\chr (k) = 2$, then $d \equiv 1- (h^*_{\Delta})^{-1} \ (\mbox{mod}\ fm_{\Delta})$ and $0 < d \leq fm_{\Delta}$.
\item In all remaining cases $d=1+2r$ where $r \equiv -(h_{\Delta})^{-1}\ (\mbox{mod}\ fm_{\Delta})$ and $0 \leq r < fm_{\Delta}$.
\end{enumerate}\end{therm}

To prove this we will first apply the results of the previous section to study how the stable Calabi-Yau property for preprojective algebras lifts to mesh algebras, and to the stable Auslander algebra $\Gamma$ in particular.  We recall from Example (2) in Section 3, that the preprojective algebra of a (simply laced) Dynkin graph $\Delta$ is the mesh algebra on the quiver $Q_{\Delta}$ obtained from $\Delta$ by replacing each edge with a pair of arrows with opposite orientations.  In fact, this translation quiver can be expressed as $\mathbb{Z}\Delta/\gen{\tau}$, and hence we can realize some of the m-fold mesh algebras as Galois covers of the preprojective algebras as outlined in Section 3.

In order to obtain information about these mesh algebras, we must consider the following ``half''-grading of $P(\Delta)$ as in \cite{Per}.  We first choose a bipartite orientation of $\Delta$, so that each vertex is either a sink or source.  We give the corresponding arrows in $Q_{\Delta}$ degree $0$, and all remaining arrows degree $1$, and consider this as a $\mathbb{Z}/\gen{m}$-grading for some positive integer $m$.  Notice that we have $\deg (\alpha) + \deg (\bar{\alpha}) = 1$ for all arrows $\alpha$ of $Q_{\Delta}$, and $\deg (\alpha \beta) = 1$ for all arrows $\alpha, \beta$ of $Q_{\Delta}$ for which $\alpha \beta \neq 0$ in $P(\Delta)$.  We shall call a vertex $v$ a $0$-source (resp. $0$-sink) if all arrows leaving (resp. entering) $v$ have degree $0$, and the notions of $1$-source and $1$-sink are defined analogously.  

The Nakayama permutation $\pi$ of $A = P(\Delta)$ is induced by the unique graph automorphism of $\Delta$ of order $2$ when $\Delta = \mathbb{A}_n\ (n \geq 2), \mathbb{D}_{2n+1}\ (n \geq 2)$ or $\mathbb{E}_6$, and is the identity otherwise (cf. Prop. 2.1 in \cite{DPA}).  We will also write $\pi$ for the corresponding quiver automorphism of $Q_{\Delta}$.  The following description of the Nakayama automorphism $\nu$ of $P(\Delta)$ is taken from Definition 4.6 of \cite{PAAK}.
\begin{eqnarray}
\nu(e_i) = \pi(e_i) \ \ (1 \leq i \leq n); \hspace{1cm} \nu(\alpha) = (-1)^{\deg(\alpha)} \pi (\alpha) \ \ \mbox{for\ all\ arrows}\ \alpha\ \mbox{in}\ Q_{\Delta}.
 \end{eqnarray}
In particular, observe that $\nu$ is a graded automorphism if and only if $\Delta$ is not $\mathbb{A}_{2n}$ for some $n$.

We now set $A = P(\Delta)$ and consider the minimal projective resolution $P_{\bullet}$ of $A$ over $A^e$, as described in (4.1-4.3).  It is not hard to see that the modules in this resolution can be graded appropriately so that it becomes a graded projective resolution (cf. Section 2 of \cite{Per}).  Furthermore, it is well known that $\Omega^3_{A^e}(A) \cong DA \cong {}_{\nu}A_1$ (see, for example, \cite{PAAK}), and it then follows that $\Omega^6_{A^e}(A) \cong {}_{\nu^2}A_1 \cong A$.  In the following proposition, we refine these to isomorphisms of graded bimodules, by studying the grading of $P_{\bullet}$.  (The reader may wish to contrast these results with those in \cite{PAAK, CYFA}, which are with respect to the natural path-length grading.)

\begin{propos} Suppose $A = P(\Delta)$ is given a half-grading as above.  Then we have the following isomorphisms of graded bimodules.
\begin{enumerate}
\item $\Omega^3_{A^e}(A) \cong DA[h_{\Delta}-1]$.
\item $\Omega^6_{A^e}(A) \cong A[h_{\Delta}]$.  
\item For any $\Delta$ other than $\mathbb{A}_{2n}$ for $n \in \mathbb{N}$, $DA[h_{\Delta}/2-1] \cong {}_{\nu} A_1$.
\end{enumerate}
\end{propos}

\noindent
{\it Proof.} We will make repeated use of the fact that all indecomposable projective $P(\Delta)$-modules have Loewy length $h_{\Delta}-1$ (See, e.g., \cite{DPA, PAAK}).

(3) Assume $\Delta$ is not $\mathbb{A}_{2n}$.  Since $\soc(A)$ is generated by pahts of length $h_{\Delta}-2$, it will be concentrated in degree $\ceil{\frac{h_{\Delta}-2}{2}} = \frac{h_{\Delta}}{2} - 1$.  Thus $DA$ is generated in degree $1-\frac{h_{\Delta}}{2}$, and we must shift it by $\frac{h_{\Delta}}{2}-1$ to have it generated in degree $0$ as is ${}_{\nu}A_1$.  (If $\Delta = \mathbb{A}_{2n}$, then $\soc(A)$ is generated by paths of odd length $2n-1$, which may have degree either $n$ or $n-1$.  Thus $DA$ is generated in degrees $-n$ and $1-n$ in this case.)

(1) Since $\Omega^3_{A^e}(A) \cong DA$ as a bimodule and it is indecomposable, it must be isomorphic to some shift of $DA$ as a graded bimodule.  (We are using the fact that when $G$ is abelian, $G$-graded bimodules correspond to $G$-graded modules over the enveloping algebra with the grading $A^e_g = \oplus_{s \in G} A_{g-s} \otimes_k A_s$.)  We thus need to determine in which degree(s) $\Omega^{3}_{A^e}(A)$ is generated.  Clearly, $P_0$ is generated in degree $0$, and $P_1$ is generated in degrees $0$ and $1$ since $\deg (x_{\alpha}) = \deg (\alpha)$ for any arrow $\alpha$.  We now claim that $P_1$ is generated in degree $1$.  The generator $e_i \otimes e_i$ of $P_2$ is mapped by $R$ to $\sigma_i \in P_1$, which is a sum of terms of the forms $\bar{\alpha} \otimes e_i \in (e_{i\alpha} \otimes e_i)A^e$ and $e_i \otimes \alpha \in (e_i \otimes e_{i\alpha})A^e$ for the different arrows $\alpha$ that end in $i$.  In either case, we see that the term has degree $\deg (\bar{\alpha}) + \deg (\alpha) = 1$.  Now $\Omega^3_{A^e}(A)$ is generated by the elements $\xi_i = \sum_{x \in e_i \mathcal{B}} x \otimes x^{*} \in P_2.$  Since $xx^{*}$, by definition, is in the socle of $e_iA$, it has degree $\frac{h_{\Delta}}{2}-1$ (when $\Delta \neq \mathbb{A}_{2n}$), or else degree $\frac{h_{\Delta} \pm 1}{2}-1$ (when $\Delta = \mathbb{A}_{2n}$).  Thus the term $x \otimes x^* \in P_2$ has degree $\frac{h_{\Delta}}{2}$ or $\frac{h_{\Delta} \pm 1}{2}$.  Comparing these degrees to those in which $DA$ is generated, yields $\Omega^3_{A^e}(A) \cong DA[h_{\Delta}-1]$.

(2) If $\Delta \neq \mathbb{A}_{2n}$, then $\Omega^3_{A^e}(A)$ is a twisted bimodule generated in degree $h_{\Delta}/2$.  Thus $\Omega^6_{A^e}$ is generated in degree $h_{\Delta}$ and hence $\Omega^6_{A^e}(A) \cong A[h_{\Delta}]$.  Now, if $\Delta = \mathbb{A}_{2n}$, we tensor $P_{\bullet}$ with a simple right $A$-module $S_i$ to get the start of a minimal projective resolution
$$\ses{\nu(S_i) \longrightarrow e_iA}{\bigoplus_{t\alpha = i} e_{i\alpha}A}{e_iA \rightarrow S_i}{}{},$$ where $\nu$ is the Nakayama functor.  Since $P_2$ is generated in degree $1$, $\nu(S_i)$ is concentrated in degree one greater than the degree of the socle of $e_iA$, so that
$$\Omega^3(S_i) = \nu(S_i)[d], \ \ \mathrm{where}\ d = \left\{ \begin{array}{cc}1 + \ceil{\frac{h_{\Delta}-2}{2}}, & i = 1\mbox{-sink} \\ 1+ \floor{\frac{h_{\Delta}-2}{2}}, & i = 1\mbox{-source}.
\end{array} \right. $$
Since  $\nu$ interchanges the $1$-sinks and $1$-sources for $P(\mathbb{A}_{2n})$, we have $\Omega^6(S_i)$ generated in degree $1+\ceil{\frac{h_{\Delta}-2}{2}} + 1 + \floor{\frac{h_{\Delta}-2}{2}} = h_{\Delta}$.  That $\Omega^6_{A^e}(A)$ is generated in degree $h_{\Delta}$ now follows from the observations in Section 2 of \cite{Per}. $\Box$ \\

We continue to let $A = P(\Delta)$ be a preprojective algebra with a $\mathbb{Z}/\gen{m}$-grading as above, so that $B = A \# k[G]^*$ is isomorphic to the mesh algebra of the translation quiver $\mathbb{Z}\Delta/ \gen{\tau^m}$.  It follows from the lemma that  $\Omega^{-3}_{A^e}(A) \cong  DA[-1]$ and $ \Omega^{-6}_{A^e}(A)  \cong  A[-h_{\Delta}]$.  Hence, we have isomorphisms of graded bimodules  
\begin{eqnarray} \Omega^{-6d}_{A^e}(A) & \cong & A[-d h_{\Delta}] \end{eqnarray} and \begin{eqnarray} \Omega^{-3-6d}_{A^e}(A) & \cong & DA[-1-dh_{\Delta}] \end{eqnarray} for all $d \geq 0$.  Lifting these isomorphisms to $B$ via the functor $F$ of the previous section, we obtain two possible cases:
\begin{eqnarray}
\Omega^{-6d}_{B^e}(B) \cong DB\ & \mbox{if\ and\ only\ if\ } & DB \cong {}_1 B_{\tau^{-dh_{\Delta}}} \cong {}_{\tau^{dh_{\Delta}}} B_1; \\ 
\Omega^{-3-6d}_{B^e}(B) \cong DB\ & \mbox{if\ and\ only\ if\ } &  dh_{\Delta} \equiv -1\ (\mbox{mod}\ m).  
\end{eqnarray}
By the remarks at the end of the previous section, we see that (6.5) can only occur if $A$ is weakly symmetric, which happens if and only if $\Delta =  \mathbb{A}_1, \mathbb{D}_{2n}, \mathbb{E}_7$ or $\mathbb{E}_8$.

\begin{therm} Let $B$ be the mesh algebra of the translation quiver $\mathbb{Z}\Delta/\gen{\tau^m}$, where $\Delta$ is a Dynkin graph.  Then there exists a $d>0$ such that $\Omega^{-3d}_{B^e}(B) \cong DB$ as bimodules if and only if $(m, h^*_{\Delta}) = 1$.  The minimal such $d$ is given as follows
\begin{enumerate}
\item If $\Delta = \mathbb{A}_1, \mathbb{D}_{2n}, \mathbb{E}_7$ or $\mathbb{E}_8$ and either $2 | m$ or $\chr (k) = 2$, then $d \equiv 1- (h^*_{\Delta})^{-1} \ (\mbox{mod}\ m)$ and $0 < d \leq m$.
\item In all remaining cases $d=1+2r$ where $r \equiv -(h_{\Delta})^{-1}\ (\mbox{mod}\ m)$ and $0 \leq r < m$.
\end{enumerate}
\end{therm}

\noindent
{\it Proof.}  We first consider the case where $P(\Delta)$ is not weakly symmetric.  As noted above, this occurs for $\Delta = \mathbb{A}_n, \mathbb{D}_{2n+1}, \mathbb{E}_6$ and $n \geq 2$.  Thus, we see from (6.6) above that $\Omega^{-3(1+2r)}_{B^e}(B) \cong DB$ if and only if $r \equiv -h_{\Delta}^{-1}\ (\mbox{mod}\ m)$.  The latter is possible if and only if $(m, h^*_{\Delta}) = 1$, and this yields the desired value of $d$.  

For the remaining $\Delta$, we break our proof into three cases based on the parity of $m$ and the characteristic of $k$.  First, assume that $\chr (k) \neq 2$.  We want to examine when (6.5) occurs; i.e., we need to determine when the action of $rh_{\Delta} \in \mathbb{Z}/\gen{m}$ induces a Nakayama automorphism of $B$.  Recall that the Nakayama automorphism $\tilde{\nu}$ for $B$ is given by $\tilde{\nu}(ap_g) = \nu(a)p_{g+h_{\Delta}/2-1}$, according to Lemmas 4.2 and 5.2, with $\nu$ as in (6.2).  It is not hard to see that $\tilde{\nu}$ differs from the automorphism induced by $h_{\Delta}/2-1 \in \mathbb{Z}/\gen{m}$ by an inner automorphism if and only if $m$ is even: the inner automorphism is conjugation by $\sum_{i \in \mathbb{Z}/\gen{m}} (-1)^i p_i$.  Furthermore, when $m$ is odd, the nontrivial outer automorphism $(h_{\Delta}/2-1)\nu^{-1}$ fixes all $e_i p_g$ and is hence not induced by any element of $\mathbb{Z}/\gen{m}$.  Hence, we see that (1) is equivalent to 
\begin{eqnarray} rh_{\Delta} & \equiv & (h_{\Delta}/2-1)\ (\mbox{mod}\ m)\end{eqnarray} 
where $m$ is even, and this yields $d = 2r$.  Under our current assumptions, $h_{\Delta}, h_{\Delta}/2 - 1$ and $m$ are all even, so this is equivalent to $rh^*_{\Delta} = (h^*_{\Delta}-1)/2 + sm/2$ for some $s \in \mathbb{Z}$.  Since $h^*_{\Delta}$ and $(h^*_{\Delta} - 1)/2$ are relatively prime, we can solve for $r$ if and only if $(h^*_{\Delta}, m/2) = 1$, which holds if and only if $(h^*_{\Delta}, m)= 1$ since $h^*_{\Delta}$ is odd.  Letting $d = 2r$ in the displayed congruence now leads to $d \equiv (h^*_{\Delta})^{-1}(h^*_{\Delta}-1) \equiv 1 - (h^*_{\Delta})^{-1}\ (\mbox{mod}\ m)$.

We continue to assume that $\chr(k) \neq 2$, and consider case (6.6).  As above, this requires that $(m, h_{\Delta}) = 1$, and since $h_{\Delta}$ is even, $m$ must be odd.  Moreover, when $m$ is odd, $(m, h_{\Delta}) =1$ if and only if $(m, h^*_{\Delta}) = 1$.  The value for $d$ in this case is obtained as above.

Finally, we assume $\chr(k)=2$.  In this case the Nakayama automorphism of $A$ is the identity, and thus $\tilde{\nu}$ is induced by $h^*_{\Delta}-1 \in \mathbb{Z}/\gen{m}$.  We thus have $\Omega^{-3}_{A^e}(A) \cong DA[-1] \cong A[-h^*_{\Delta}]$, from which we obtain $\Omega^{-3d}_{B^e}(B) \cong {}_1 B_{\tau^{-dh^*_{\Delta}}} \cong {}_{\tau^{dh^*_{\Delta}}}B_1$.  We now see that the latter is isomorphic to $DB \cong {}_{\tau^{h^*_{\Delta}-1}}B_1$ if and only if $dh^*_{\Delta} \equiv h^*_{\Delta}-1 \ (\mbox{mod}\ m)$.  Clearly, this is possible if and only if $(h^*_{\Delta}, m) = 1$, and in this case we again obtain $d \equiv 1-(h^*_{\Delta})^{-1}\ (\mbox{mod}\ m)$.  $\Box$\\

Theorem 6.1 now follows immediately from the above theorem, using Proposition 2.1 and the description of the mesh algebras given at the beginning of the section.

\section{Generalized preprojective algebras and the torsion order 2 case}
\setcounter{equation}{0}

We now move on to study the stable Calabi-Yau properties of the representation-finite self-injective algebras having torsion order $t=2$.  These algebras fall into three classes corresponding to their Dynkin type.  More specifically, the possible types here are $(\mathbb{A}_{2n+1}, s, 2),\  (\mathbb{D}_n, s, 2)$ and $(\mathbb{E}_6, s, 2)$ with $s \in \mathbb{N}$.  Unfortunately, the lifting methods used above are not as well-suited for this case.  Already, in \cite{Per}, we found that we could not find exact values for the periods of these algebras, due to the possibility that some syzygy $\Omega_{A^e}(A)$ could be a nontrivially twisted bimodule which lifts to the regular bimodule over a double cover of the algebra.   In general, there do exist nontrivial outer automorphisms that lift to inner autmorphisms in this way, and hence could provide such a twisting, but we know of no elementary means of detecting them in syzygies.

For now, we avoid this problem altogether by lifting the minimal resolutions of generalized preprojective algebras to obtain information about the relevant stable Auslander algebras.  The drawback of this approach is that not all the stable Auslander algebras that occur are coverings of generalized preprojective algebras.  We thus will only be able to derive some partial results for the time being.  In particular, for the cases $(\mathbb{A}_{2n+1},s,2)$ and $(\mathbb{D}_n, s, 2)$, we obtain precise results concerning stable Calabi-Yau dimensions only when $s$ is odd.  We also get some new (although still incomplete) information about the periods of algebras of types $(\mathbb{D}_n, s, 2)$ and $(\mathbb{E}_6, s, 2)$.

We let $\Gamma$ be the stable Auslander algebra of $\Lambda$, which is isomorphic to the mesh algebra of the translation quiver $\mathbb{Z}\Delta/\gen{\rho \tau^{sm_{\Delta}}}$, where $\rho^2 = 1$.  We let $B$ be the mesh algebra of $\mathbb{Z}\Delta/\gen{\tau^{2sm_{\Delta}}}$, which is a $2$-fold Galois cover of $\Gamma$, with Galois group generated by $x = \rho \tau^{sm_{\Delta}}$.  We thus consider $\Gamma$ with the corresponding $\mathbb{Z}/\gen{2}$-grading.  As in the last section, we also consider $B$ as a Galois cover of the preprojective algebra $A = P(\Delta)$.  As in the last section, we give $\Delta$ a bipartite orientation and assign the associated ``half''-grading to $A$.  This gives the following schemata of coverings between translation quivers and their mesh algebras.

$$\begin{array}{ccccc} \mathbb{Z}\Delta/\gen{\rho \tau^{m_{\Delta}s}} & \longleftarrow & \mathbb{Z}\Delta/\gen{\tau^{2m_{\Delta}s}} & \longrightarrow &\mathbb{Z}\Delta/\gen{\tau} \\ \\ \Gamma & & B & & A = P(\Delta) \end{array}$$

Now suppose $\Omega^{-3d}_{\Gamma^e}(\Gamma) \cong D\Gamma$ for some $d>0$.  As graded bimodules, we either have $D\Gamma$ or $D\Gamma[x]$.  Hence, lifting to $B$ we obtain $\Omega^{-3d}_{B^e}(B) \cong DB$ or $\Omega^{-3d}_{B^e}(B) \cong DB_x$.  We begin each of the following cases by applying our results from the previous section to investigate for which values of $d$ and $s$ these isomorphisms hold.  However, such isomorphisms alone do not guarantee that $\Omega^{-3d}_{\Gamma^e}(\Gamma) \cong D\Gamma$, as it is not hard to find outer automorphisms $\sigma$ of $\Gamma$ such that $F(D\Gamma_{\sigma}) \cong F(D\Gamma)$ as bimodules.

\subsection{Type $(\mathbb{A}_{2n+1}, s, 2)$}

We now assume that $\Lambda$ is a representation-finite self-injective algebra of type $(\mathbb{A}_{2n+1},s,2)$ with $n, s \geq 1$.  The AR-quiver of $\Lambda$ has the form $\mathbb{ZA}_{2n+1}/\gen{\tau^{(2n+1)s}\rho}$ where $\rho$ is the automorphism given by reflection in the horizontal central line of $\mathbb{ZA}_{2n+1}$.  We set $m = s(2n+1)$ and $x = \rho \tau^m$.  Proposition 6.2 implies that $\Omega^{-3}_{B^e}(B) \cong {}_1 DB_{\tau^{-1}}$ and $\Omega^{-6}_{B^e}(B) \cong {}_1 B_{\tau^{-(2n+2)}}$.  Hence, for $r \geq 0$ \begin{eqnarray} \Omega^{-6r}_{B^e}(B) \cong {}_1 B_{\tau^{-r(2n+2)}}\ \ \mbox{and}\ \  \Omega^{-3-6r}_{B^e}(B) \cong {}_1 DB_{\tau^{-1-r(2n+2)}}.\end{eqnarray}  

\begin{lemma} If $\Lambda$ (as above) is stably $d$-Calabi-Yau, then $s \equiv n\ (\mbox{mod}\ 2)$ and $d$ is even.
\end{lemma}

\noindent
{\it Proof.}  We suppose that $\Omega^{-3d}_{\Gamma^e}(\Gamma) \cong D\Gamma$.  Notice that as a graded bimodule, $\Omega^{-3d}(\Gamma)$ cannot be isomorphic to $D\Gamma$, for such an isomorphism would lift to an isomorphism $\Omega^{-3d}(B) \cong DB$ contradicting the requirement of Theorem 6.3 that $(2m,h_{\Delta}^*)=1$.  Thus we must have $\Omega^{-3d}(\Gamma) \cong D\Gamma[x]$, which lifts to 
$$\Omega^{-3d}(B) \cong DB_x \cong DB_{\rho\tau^m} \cong {}_{\rho\tau^n} B_{\rho\tau^m} \cong {}_1 B_{\tau^{m-n}}.$$  We have used that the Nakayama automorphism of $B$ is $\nu \tau^n$ by Lemma 5.2 and Proposition 6.2, where $\nu$ represents the extension of the Nakayama automorphism of $A$.  Furthermore, although $\nu$ multiplies some arrows by $-1$, it is not hard to see that it differs from $\rho$ by an inner automorphism since the quiver of $B$ is a wreath with even circumference.

As shown above in (7.1), a twist of $B$ by a power of $\tau$ can only occur in a $6r^{th}$-cosyzygy.  Hence $d$ must be even, and we must have 
\begin{eqnarray} (-n-1)d \equiv m-n\ (\mbox{mod}\ 2m), \end{eqnarray} which forces $m-n$ to be even. Thus $m =(2n+1)s \equiv s \equiv n\ (\mbox{mod}\ 2)$.  $\Box$ \\

We now assume that $s$ and $n$ are odd.  In this case the stable AR-quiver of $\Lambda$ is an $m$-fold covering of the stable translation quiver $Q_{\mathbb{B}_{n+1}}$, whose mesh algebra $A$ is the generalized preprojective algebra $P(\mathbb{B}_{n+1})$.  As for the preprojective algebras, we can give $P(\mathbb{B}_{n+1})$ a half-grading corresponding to this covering by assigning degrees of $0$ and $1$ to the arrows as indicated below.  The figure on the right shows a piece of $\mathbb{ZA}_{2n+1}$, with the vertex labels illustrating the covering morphism and the solid (resp. dashed) arrows covering the arrows of degree $0$ (resp. $1$).

$$\begin{array}{ccc} & & \vdots \\ 
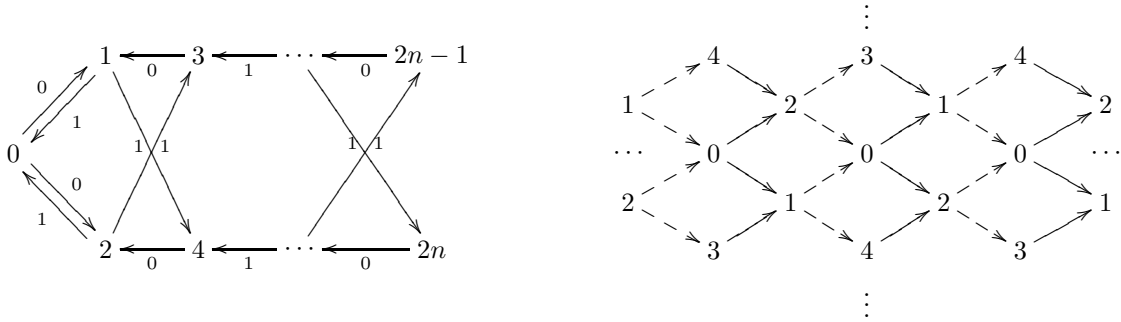
\xymatrix{& 1 \ar[ddr]^1  \ar[dl]<0.5ex>^1 & 3 \ar[l]^0  & \cdots \ar[l]^1 \ar[ddr]^1 & 2n-1 \ar[l]^0  \\ 0 \ar[ur]<0.5ex>^0 \ar[dr]<0.5ex>^0 \\ & 2 \ar[ul]<0.5ex>^1 \ar[uur]^1 & 4 \ar[l]^0 & \cdots \ar[l]^1  \ar[uur]^1 & 2n \ar[l]^0} & \hspace{1cm} &  
\xymatrixrowsep{0.5pc} \xymatrixcolsep{1.5pc} \xymatrix{ & 4 \ar[dr] & & 3 \ar[dr] & & 4 \ar[dr] & \\ 1  \ar@{-->}[ur] \ar@{-->}[dr] & & 2 \ar@{-->}[ur] \ar@{-->}[dr] & & 1 \ar@{-->}[ur] \ar@{-->}[dr] & & 2  \\ \cdots & 0 \ar[ur] \ar[dr] & & 0 \ar[ur] \ar[dr]  & & 0 \ar[ur] \ar[dr] & \cdots \\ 2 \ar@{-->}[ur] \ar@{-->}[dr]  & & 1 \ar@{-->}[ur] \ar@{-->}[dr] & & 2 \ar@{-->}[ur] \ar@{-->}[dr] & & 1 \\ & 3 \ar[ur]  & & 4 \ar[ur]  & & 3 \ar[ur]  &  } 
 \\  & & \vdots \end{array}$$

Our task is simplified somewhat by the fact that, for $n$ odd, $P(\mathbb{B}_{n+1})$ appears as the stable Auslander algebra of the category $\mathcal{C}(R)$ of Cohen-Macaulay modules over the simple plane curve singularity $R = k[[x,y]]/(x^2y + y^{n+1})$ of type $\mathbb{D}_{n+2}$.  In the stable category $\und{\mathcal{C}}(R)$, one knows that $\Omega \cong \tau$ and $\Omega^2 \cong Id$.  Furthermore, $\Omega \tau \cong Id$ is a Serre functor, and we see that $\und{\mathcal{C}}(R)$ is a $2$-Calabi-Yau triangulated category (see \cite{Yosh}).  Hence, Proposition 2.1 implies that $\Omega^{-6}_{A^e}(A) \cong A \cong DA$.

  We must now analyze the resolution described in Section 4 in order to determine the degrees in which the syzygies are generated.  Note that since $\mathbb{ZA}_{2n+1}$ is the universal cover of $Q_{\mathbb{B}_{n+1}}$, the indecomposable projective $A$-modules will have the same Loewy length as the indecomposable projective $P(\mathbb{A}_{2n+1})$-modules, which is $h_{\mathbb{A}_{2n+1}}-1 = 2n+1$.  Thus $DA$ is generated in degree $-n$ and we have a graded isomorphism $DA[n] \cong A$.  As demonstrated for the preprojective algebras in the proof of Proposition 6.2, $P_2$ is generated in degree $1$.  Furthermore, the generators $\xi_i$ of $L = \Omega^3(A)$ are each homogeneous of degree $1+n$, so $\Omega^6_{A^e}(A)$ will be generated in degree $2n+2 = h_{\mathbb{A}_{2n+1}}$.  This yields $\Omega^{-6r}_{A^e}(A) \cong A[-(2n+2)r] \cong DA[n-(2n+2)r]$, and applying the lifting functor $F$ gives us $\Omega^{-6r}_{\Gamma^e}(\Gamma) \cong {}_1 D\Gamma_{\tau^{n-(2n+2)r}}$.  Thus $\stmod \Lambda$ is $2r$-Calabi-Yau if and only if $(2n+2)r \equiv n\ (\mbox{mod}\ s(2n+1))$.  This is possible if and only if $(2n+2, s) | n$, and since $(2n+2,n)=1$ (as $n$ is odd), this condition is equivalent to $(2n+2,s) =1$ or even $(n+1, s) = 1$ (as $s$ is odd).
  
  We summarize what we have found in the following proposition.  The very last assertion follows from the congruence (7.2).  The uncertainty in the case where $n$ is even is resolved in Proposition 9.6.
  
\begin{propos} Let $\Lambda$ be a representation-finite self-injective algebra of type $(\mathbb{A}_{2n+1}, s, 2)$.
\begin{enumerate}
\item If $n$ is odd, then $\stmod \Lambda$ is $d$-Calabi-Yau if and only if $(n+1,s)=1$.  In this case $d = 2r$ with $r \equiv n(2n+2)^{-1}\ (\mbox{mod}\ s(2n+1))$ and $0 < r < s(2n+1)$.
\item If $n$ is even and $\stmod \Lambda$ is $d$-Calabi-Yau, then $s$ is even and $(n+1,s)=1$.  $\Box$
\end{enumerate}
\end{propos}  

\subsection{Type $(\mathbb{D}_n, s, 2)$}

Now suppose $\Lambda$ has type $(\mathbb{D}_n, s, 2)$.  The stable Auslander algebra $\Gamma$ of $\Lambda$ coincides with the mesh algebra of $\mathbb{ZD}_n/\gen{\rho \tau^{(2n-3)s}}$, where $\rho$ is the automorphism induced by the order 2 graph automorphism of $\mathbb{D}_n$.  As above, we let $B$ be the mesh algebra of $\mathbb{ZD}_n/\gen{\tau^{2s(2n-3)}}$ which is a degree-$2$ cover of $\Gamma$ and a degree-$2s(2n-3)$ cover of $A = P(\mathbb{D}_n)$.  By Proposition 6.2 we have $\Omega^{-3}_{A^e}(A) \cong DA[-1]$ and $\Omega^{-6}_{A^e}(A) \cong A[2-2n]$.  Lifting to $B$, we have 
\begin{eqnarray} \Omega^{-6r}(B) \cong {}_1 B_{\tau^{(2-2n)r}} &  \mbox{and} &  \Omega^{-3-6r}(B) \cong DB_{\tau^{-1 +r(2-2n)}} \cong {}_{\nu \tau^{n-2}} B_{\tau^{-1+r(2-2n)}} \cong {}_1 B_{\nu \tau^{1-n+r(2-2n)}}. \end{eqnarray}
  Here, $\nu$ is induced by the Nakayama automorphism of $A$, which is the identity if $n$ is even, and corresponds to the order $2$ automorphism $\rho$ of $\mathbb{D}_n$ when $n$ is odd.  We have also used that the Nakayama automorphism of $B$ is $\tilde{\nu} = \nu\tau^{n-2}$ according to Lemma 5.2 and Proposition 6.2.  While $\nu$ multiplies the degree-one arrows by $-1$, since the quiver of $B$ is a wreath of even circumference we can compose $\tilde{\nu}$ with an inner automorphism so that it takes arrows to arrows.  We may thus assume that $\tilde{\nu} = \tau^{n-2}$ or $\rho\tau^{n-2}$ when $n$ is even or odd, respectively.

Now suppose that $\Omega^{-3d}(\Gamma) \cong D\Gamma$.  To lift this to an isomorphism of $B$-bimodules we must consider two cases.  If $\Omega^{-3d}(\Gamma) \cong D\Gamma$ as $\mathbb{Z}/\gen{2}$-graded bimodules, then $\Omega^{-3d}(B) \cong DB$; while if $\Omega^{-3d}(\Gamma) \cong D\Gamma[x]$ as $\mathbb{Z}/\gen{2}$-graded bimodules, then $\Omega^{-3d}(B) \cong {}_1 DB_{\rho \tau^{s(2n-3)}}$, where $x$ is the generator of $\mathbb{Z}/\gen{2}$ and corresponds to the automorphism $\rho \tau^{s(2n-3)}$ of $\mathbb{ZD}_n/\gen{\tau^{2s(2n-3)}}$.  We first consider the case $\Omega^{-3d}(B) \cong DB$.  By Theorem 6.3, we would have $(2s(2n-3), 2n-2) = 1$ when $n$ is odd -- which is clearly false -- or else $(2s(2n-3), n-1) = 1$ when $n$ is even.  Note that the latter is then equivalent to $(s,n-1) = 1$.  Theorem 6.3(1) also yields $d \equiv 1 - (n-1)^{-1}\ (\mbox{mod}\ 2s(2n-3))$, which forces $d$ to be even.

In the remaining case, $\Omega^{-3d}(B) \cong {}_1 DB_{\rho \tau^{s(2n-3)}} \cong {}_{\tilde{\nu}}B_{\rho \tau^{s(2n-3)}}$.  If $n$ is even, this module is isomorphic to ${}_1 B_{\rho \tau^{s(2n-3)-n+2}}$.  However, this cannot occur since (7.3) shows that $\Omega^{-6r}(B) \cong {}_1B_{\tau^{r(2-2n)}}$ and $\Omega^{-3-6r}(B) \cong {}_1 B_{\tau^{1-n+r(2-2n)}}$ when $n$ is even, and $\rho$ is not an inner automorphism.  On the other hand, if $n$ is odd, we get $\Omega^{-3d}(B) \cong {}_1B_{\tau^{s(2n-3)-n+2}}$.  It follows from (7.3) that $d$ is even and $s(2n-3)-n+2 \equiv (1-n)d\ (\mbox{mod}\ 2s(2n-3))$.  In particular, we see that $s(2n-3) - n + 2$ is even.  Consequently, if $n$ is odd, then $s$ too must be odd.

We now assume that $s$ is odd, but do not place any restrictions on $n$.  In this case the stable AR-quiver of $\Lambda$, which has the form $\mathbb{ZD}_n/\gen{\rho \tau^{s(2n-3)}}$, is an $s(2n-3)$-fold cover of the translation quiver $Q_{\mathbb{C}_{n-1}}$, which is obtained by factoring out the action of the cyclic group of automorphisms genereated by $\rho \tau$.  The stable Auslander algebra $\Gamma$ is thus a Galois cover of the generalized preprojective algebra $P(\mathbb{C}_{n-1})$, and the corresponding grading on the latter algebra is indicated below.  On the right is  a piece of the AR-quiver of $\Lambda$ with the vertex labels showing the covering morphism.

$$\begin{array}{ccc}  
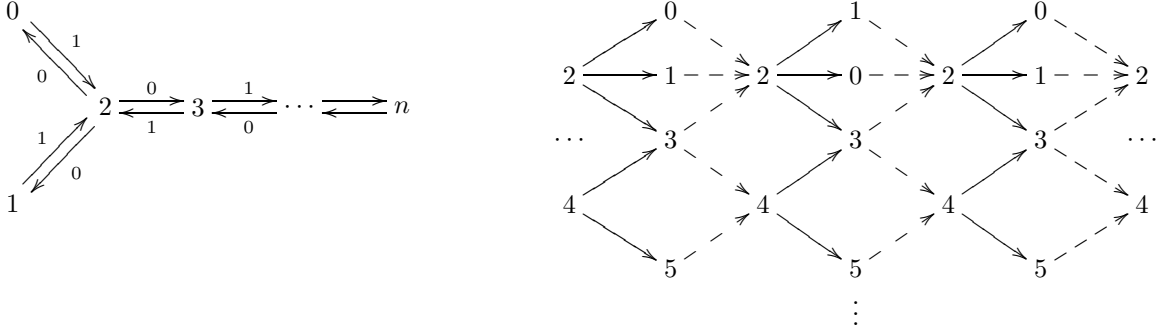
\xymatrix{0 \ar[dr]<0.5ex>^1  \\ & 2 \ar[ul]<0.5ex>^0 \ar[r]<0.5ex>^0 \ar[dl]<0.5ex>^0 & 3 \ar[r]<0.5ex>^1 \ar[l]<0.5ex>^1 & \cdots \ar[l]<0.5ex>^0 \ar[r]<0.5ex> &  n \ar[l]<0.5ex> \\ 1 \ar[ur]<0.5ex>^1} 
& \hspace{1cm} &  
\xymatrixrowsep{1.0pc} \xymatrixcolsep{2.0pc} \xymatrix{ & 0 \ar@{-->}[dr] & & 1 \ar@{-->}[dr] & & 0 \ar@{-->}[dr] & \\ 2  \ar[ur] \ar[dr] \ar[r] & 1 \ar@{-->}[r] & 2 \ar[ur] \ar[dr] \ar[r] & 0 \ar@{-->}[r] & 2 \ar[ur] \ar[dr] \ar[r] & 1 \ar@{-->}[r] & 2  \\ \cdots & 3 \ar@{-->}[ur] \ar@{-->}[dr] & & 3 \ar@{-->}[ur] \ar@{-->}[dr]  & & 3 \ar@{-->}[ur] \ar@{-->}[dr] & \cdots \\ 4 \ar[ur] \ar[dr]  & & 4 \ar[ur] \ar[dr] & & 4 \ar[ur] \ar[dr] & & 4 \\ & 5 \ar@{-->}[ur]  & & 5 \ar@{-->}[ur]  & & 5 \ar@{-->}[ur]  &  } 
 \\  & & \vdots \end{array}$$

As in our previous analysis, we use the fact that $Q_{\mathbb{C}_{n-1}}$ is the stable AR-quiver of the category $\mathcal{C}(R)$ of Cohen-Macaulay modules over the simple plane curve singularity $R = k[[x,y]]/(x^2 + y^{2n-2})$ of type $\mathbb{A}_{2n-3}$.  Thus, setting $A := P(\mathbb{C}_{n-1})$, we have as above $\Omega^{-6}_{A^e}(A) \cong DA \cong A$.  To calculate the grading on the minimal resolution $P_{\bullet}$ of $A$, we observe that the indecomposable projective $A$-modules have Loewy length $h_{\mathbb{D}_n} -1 = 2n-3$.  It follows that $DA$ is generated in degree $-(n-2)$ and hence $DA \cong A[2-n]$.  Next, observe that as for the preprojective algebras, $P_1$ is generated in degree $1$.  To find the degree of the generators $\xi_i$ of $\Omega^3_{A^e}(A)$, notice that each term of the form $\tau(x) \otimes x^*$ as in (\ref{eq:xidef}) will have degree equal to $(|x| + |x^*|)/2 = n-2$, since the numerator is always one less than the Loewy length $2n-3$ of an indecomposable projective.  Thus $\Omega^3_{A^e}(A)$ is generated in degree $n-1$, and we must then have $\Omega^6_{A^e}(A) \cong A[2n-2]$.  Now we will have $\Omega^{-6r}_{A^e}(A) \cong DA$ as graded bimodules, or equivalently $\Omega^{-6r}_{\Gamma^e}(\Gamma) \cong D\Gamma$, if and only if $(2-2n)r \equiv 2-n\ (\mbox{mod}\ s(2n-3))$.  Thus we obtain the following result.

\begin{propos} Let $\Lambda$ be a representation-finite self-injective algebra of type $(\mathbb{D}_n, s, 2)$.
\begin{enumerate}
\item If $s$ is odd, then $\stmod \Lambda$ is $d$-Calabi-Yau if and only if $(n-1,s)=1$.  In this case the minimal such $d$ is  $2r$ with $r \equiv (n-2)(2n-2)^{-1}\ (\mbox{mod}\ s(2n-3))$ and $0 < r < s(2n-3)$.
\item If $s$ is even and $\stmod \Lambda$ is $d$-Calabi-Yau, then $(n-1,s)=1$.
\item If $s>1$ is odd, then the period of $\Lambda$ is $\frac{2s(2n-3)}{(s,n-1)}$. $\Box$
\end{enumerate}
\end{propos}

For (3), we have used the isomorphism $\Omega^6_{\Gamma^e}(\Gamma) \cong {}_1 \Gamma_{\tau^{2n-2}}$ to get $\Omega^{6s(2n-3)/(s,2n-2)}_{\Gamma^e}(\Gamma) \cong {}_1 \Gamma_1$.  Thus $\Gamma$ has period dividing $6s(2n-3)/(s,n-1)$, and hence $\Lambda$ has period dividing $2s(2n-3)/(s,n-1)$ by Theorem 4.2(2) of \cite{Per} (this requires $s>1$ so that $\Lambda$ is Schurian).  But this coincides with the lower bound for the period of $\Lambda$ that was computed in \cite{Per}.  In Proposition 9.7, we supplement parts (2) and (3) with additional information in the case where $s$ is even and $n$ is odd.

\subsection{Type $(\mathbb{E}_6, s, 2)$}

We now assume $\Lambda$ is a representation-finite self-injective algebra of type $(\mathbb{E}_6, s, 2)$.  The stable AR-quiver of $\Lambda$ then has the form $\mathbb{ZE}_6/\gen{\rho \tau^{11s}}$ where $\rho$ is the automorphism of $\mathbb{ZE}_6$ induced by the order-$2$ graph automorphism of $\mathbb{E}_6$.  The stable Auslander algebra $\Gamma$ of $\Lambda$ is the mesh algebra of this translation quiver.  As above, we begin by looking at the double cover $B$ of $\Gamma$, given by the mesh algebra of $\mathbb{ZE}_6/\gen{\tau^{22s}}$, and we let $A = P(\mathbb{E}_6)$ for the time being.

Assume that $\Omega^{-3d}_{\Gamma^e}(\Gamma) \cong D\Gamma$.  If this is an isomorphism of $\mathbb{Z}/2$-graded bimodules, then it lifts to an isomorphism $\Omega^{-3d}_{B^e}(B) \cong DB$, which leads to a contradiction $(22s, 12)=1$ by Theorem 6.3.  Thus, we must have instead $\Omega^{-3d}_{B^e}(B) \cong {}_1 DB_{\rho \tau^{11s}}$.  Now, $DB \cong {}_{\tilde{\nu}}B_1 \cong {}_{\rho \tau^5} B_1$ by Lemma 5.2 and Proposition 6.2.  Thus $\Omega^{-3d}_{B^e}(B) \cong {}_1 B_{\tau^{11s-5}}$.  Over $A$, we have $\Omega^{-6r}_{A^e}(A) \cong A[-12r]$ and $\Omega^{-3-6r}_{A^e}(A) \cong DA[-1-12r]$.  Only the first of these isomorphisms will lift to one of the desired form over $B$, and we see that $d=2r$ will be even, and we need $-12r \equiv 11s-5\ (\mbox{mod}\ 22s)$.  In particular, we see that $s$ must be odd.

Since we now know that $s$ is odd, we can regard $\Gamma$ as a Galois cover of the generalized preprojective algebra $A = P(\mathbb{F}_4)$.  The relevant half-grading of $A$ is given below, and the figure on the right shows a piece of the AR-quiver of $\Lambda$.  As above, the vertex labels on the right define the covering morphism.

$$\begin{array}{ccc}  
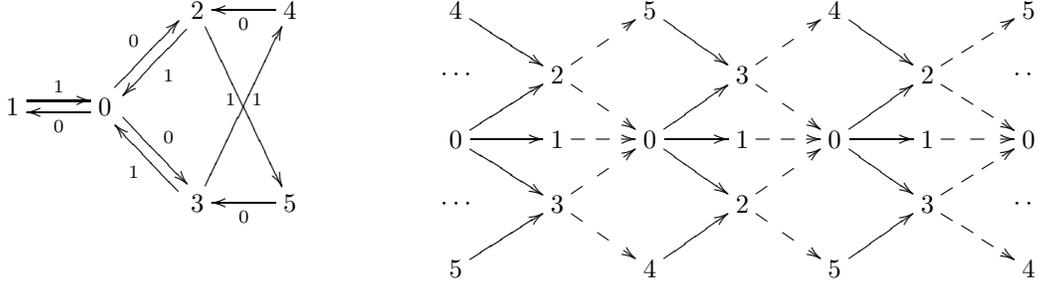
\xymatrix{ & & 2 \ar[dl]<0.5ex>^1  \ar[ddr]^1 & 4 \ar[l]^0  \\ 1 \ar[r]<0.5ex>^1 & 0 \ar[l]<0.5ex>^0 \ar[ur]<0.5ex>^0 \ar[dr]<0.5ex>^0 \\ & & 3 \ar[ul]<0.5ex>^1 \ar[uur]^1 & 5 \ar[l]^0} 
& \hspace{1cm} &  
\xymatrixrowsep{1.0pc} \xymatrixcolsep{2.0pc} \xymatrix{ 4 \ar[dr] & & 5 \ar[dr] & & 4 \ar[dr] & & 5 \\ \cdots & 2 \ar@{-->}[dr] \ar@{-->}[ur] & & 3 \ar@{-->}[ur] \ar@{-->}[dr] & & 2 \ar@{-->}[ur] \ar@{-->}[dr] & \cdots \\ 0  \ar[ur] \ar[dr] \ar[r] & 1 \ar@{-->}[r] & 0 \ar[ur] \ar[dr] \ar[r] & 1 \ar@{-->}[r] & 0 \ar[ur] \ar[dr] \ar[r] & 1 \ar@{-->}[r] & 0  \\ \cdots & 3 \ar@{-->}[ur] \ar@{-->}[dr] & & 2 \ar@{-->}[ur] \ar@{-->}[dr]  & & 3 \ar@{-->}[ur] \ar@{-->}[dr] & \cdots \\ 5 \ar[ur]  & & 4 \ar[ur]  & & 5 \ar[ur]  & & 4 }  \end{array}$$

Since $A$ is the stable Auslander algebra of the category of $\mathcal{C}(R)$ of Cohen-Macaulay modules over the plane curve singularity $R=k[[x,y]]/(x^3+y^4)$ of type $\mathbb{E}_6$, we again have $\Omega^6_{A^e}(A) \cong A \cong DA$.  As above, one can check using the minimal resolution of $A$ that $\Omega^3_{A^e}(A)$ is generated in degree $1 + (h_{\mathbb{E}_6}/2-1) = 6$ and hence $\Omega^{-6}_{A^e}(A) \cong A[-12]$.  One also checks as before that $DA \cong A[-5]$.  It now follows that, in order for $\Omega^{-6r}_{\Gamma^e}(\Gamma) \cong D\Gamma$, we need $12r \equiv 5\ (\mbox{mod}\ 11s)$.  Summarizing, we have the following.

\begin{propos} Let $\Lambda$ be a representation-finite self-injective algebra of type $(\mathbb{E}_6, s, 2)$.
\begin{enumerate}
\item $\stmod \Lambda$ is $d$-Calabi-Yau for some $d$ if and only if $(s,6)=1$.  In this case the minimal $d$ is given by $d = 2r$, where $r \equiv 5\cdot 12^{-1}\ (\mbox{mod}\ 11s)$ and $0 < r < 11s$.
\item If $s>1$ is odd, then the period of $\Lambda$ is $22s/(s,6)$. $\Box$
\end{enumerate}
\end{propos}

As in 7.2, the period of $\Lambda$ in part (2) is obtained by finding the minimal $p>0$ such that $\Omega^{3p}_{A^e}(A) \cong A$ as $\mathbb{Z}/\gen{11s}$-graded bimodules.  We compute the periods when $s$ is even in Proposition 9.8.

\section{Some results on $\mathbb{G}_2$-mesh algebras}
\setcounter{equation}{0}

It is already known that the representation-finite self-injective algebras of type $(\mathbb{D}_4, s, 3)$ are not stably Calabi-Yau for any $s$ \cite{BiaSko}.   We do not repeat the proof here: one simply checks that no power of the syzygy operator induces the Nakayama permutation on the simple modules.  However the $m$-fold mesh algebras of Dynkin type $\mathbb{G}_2$ are still of interest.  In particular, knowledge of their homology would help determine the periods of the self-injective algebras of type $(\mathbb{D}_4, s, 3)$ -- a question which was left unresolved in \cite{Per}.  Thus we devote this section to investigating the minimal graded resolution of the generalized preprojective algebra $P(\mathbb{G}_2)$, and to describing what information we can extract about the $m$-fold mesh algebras of type $\mathbb{G}_2$.  We point out that the translation quiver $Q_{\mathbb{G}_2}$ does not arise as the stable AR-quiver of the category of Cohen-Macaulay modules over a simple curve singularity, 
 so we must start from scratch.

We set $A = P(\mathbb{G}_2)$, which is the mesh algebra of the translation quiver below.  

$$
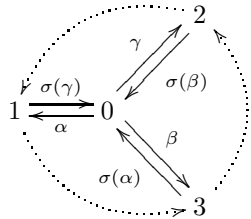
\xymatrix{  & & 2 \ar[dl]<0.5ex>^{\sigma(\beta)} \ar@{<.}@/^1.5pc/[dd] \\ 1 \ar[r]<0.5ex>^{\sigma(\gamma)} \ar@{<.}@/^1.5pc/[urr] & 0 \ar[l]<0.5ex>^{\alpha} \ar[ur]<0.5ex>^{\gamma} \ar[dr]<0.5ex>^{\beta} \\ & & 3 \ar[ul]<0.5ex>^{\sigma(\alpha)} \ar@{<.}@/^1.5pc/[ull]}$$

Notice that the Nakayama permutation $\pi$ of $A$ agrees with $\tau$ on vertices.  Each indecomposable projective $A$-module has Loewy length $h_{\mathbb{D}_4}-1 = 5$.  To lift our results to $m$-fold mesh algebras of type $\mathbb{G}_2$, we will give $A$ a $G=\mathbb{Z}/\gen{m}$ grading by assigning degree $0$ to the arrows $\alpha, \beta, \gamma$ and degree $1$ to the remaining arrows.  If $3 \nmid m$, we have $A \# k[G]^* \cong k(\mathbb{ZD}_4/\gen{\rho \tau^m})$ where $\rho$ is induced by an order $3$ automorhpism of $\mathbb{D}_4$.  If $3 \mid m$, then $A \# k[G]^* \cong k(\mathbb{ZD}_4/\gen{\tau^m})$.

\begin{propos} $\Omega^3_{A^e}(A) \cong {}_1 A_{\mu}$ where $\mu$ is the automorphism of $A$ that sends $\alpha \mapsto - \alpha, \beta \mapsto -\beta$ and $\gamma \mapsto -\gamma$, and fixes the remaining arrows and vertices of $Q_{\mathbb{G}_2}$.  Moreover, with the grading described above, $\Omega^6_{A^e}(A) \cong A[6]$; while $\Omega^3_{A^e}(A) \cong {}_1 A_{\mu}[3]$, which is isomorphic to $A[3]$ if and only if $\chr (k) = 2$.
\end{propos}

\noindent
{\it Proof.}  To find $\Omega^3_{A^e}(A)$ using the resolution from Section 4, we must first calculate its generators $\zeta_i$ for $i = 0, 1, 2, 3$.  Since these are defined in terms of a pair of dual bases with respect to the bilinear form $(-,-)$, we specify such a pair in Table 8.1.

For $\xi = \xi_0 + \cdots + \xi_3$, one can now compute $\xi \cdot a = \mu(a) \xi$ for each arrow $a$ in $Q_{\mathbb{G}_2}$, where $\mu$ is as described above.  Thus $1 \mapsto \xi$ defines an isomorphism ${}_1 A_{\mu} \cong \Omega^3_{A^e}(A)$.  Clearly $\mu$ is a nontrivial outer automorphism unless $\chr (k)=2$.  Furthermore, as in the proof of Proposition 6.2, we see that the $P_1$ term of the minimal resolution of $A$ is generated in degree $1$, and each $\xi_i$ has degree $2$.  Thus $\Omega^3_{A^e}(A)$ is generated in degree $3$. $\Box$ \\

\setcounter{table}{0}
\renewcommand{\thetable}{\thesection .\arabic{table}}
\renewcommand{\arraystretch}{1.2}
\begin{table}
\begin{tabular}{|l@{\hspace{-2mm}}r|l@{\hspace{-2mm}}r|l@{\hspace{-2mm}}r|l@{\hspace{-2mm}}r|} \hline
$ e_0\mathcal{B}$ & $\mathcal{B}^* e_0$ & $e_1 \mathcal{B}$ & $\mathcal{B}^*e_3$ & $e_2 \mathcal{B}$ & $\mathcal{B}^*e_1$ & $e_3 \mathcal{B}$ & $\mathcal{B}^* e_2 $ \\
  \hline $ e_0$ & $\alpha\sigma(\gamma)\alpha \sigma(\gamma)$ & $e_1$ & $\sigma(\gamma)\alpha \sigma(\gamma) \beta$ & $e_2$ & $\sigma(\beta)\gamma \sigma(\beta)\alpha$ & $e_3$ & $\sigma(\alpha)\beta \sigma(\alpha)\gamma$ \\
  $\alpha$ & $\sigma(\gamma)\alpha \sigma(\gamma)$ & $\sigma(\gamma)$ & $\alpha \sigma(\gamma) \beta$ & $\sigma(\beta)$ & $\gamma \sigma(\beta) \alpha$ & $\sigma(\alpha)$ & $\beta \sigma(\alpha)\gamma$ \\
  $\beta$ & $\sigma(\alpha)\beta \sigma(\alpha)$ & $\sigma(\gamma)\alpha $ & $\sigma(\gamma) \beta$ & $\sigma(\beta)\gamma$ &  $\sigma(\beta) \alpha$ & $\sigma(\alpha)\beta$ & $\sigma(\alpha)\gamma$\\
  $\gamma$ & $\sigma(\beta)\gamma \sigma(\beta)$ & $\sigma(\gamma)\beta $ & $-\sigma(\alpha) \beta$ & $\sigma(\beta)\alpha$ &  $-\sigma(\gamma) \alpha$ & $\sigma(\alpha)\gamma$ & $-\sigma(\beta)\gamma$\\
  $\alpha \sigma(\gamma)$ & $\alpha \sigma(\gamma)$ & $\sigma(\gamma)\alpha \sigma(\gamma)$ & $\beta$ & $\sigma(\beta)\gamma \sigma(\beta)$ & $\alpha$ & $\sigma(\alpha)\beta\sigma(\alpha)$ & $\gamma$\\
  $\beta \sigma(\alpha)$ & $-\gamma \sigma(\beta)$ & $\sigma(\gamma)\alpha \sigma(\gamma) \beta$ & $e_3$ & $\sigma(\beta)\gamma \sigma(\beta) \alpha$ & $e_1$ & $\sigma(\alpha)\beta \sigma(\alpha)\gamma$ & $e_2$\\
  $\gamma \sigma(\beta) \gamma $ & $\sigma(\beta)$ & & & & & & \\
$\alpha \sigma(\gamma)\alpha$ & $\sigma(\gamma) $ & & & & & &\\
$\beta \sigma(\alpha) \beta$ & $\sigma(\alpha) $ & & & & & &\\
$\alpha \sigma(\gamma)\alpha \sigma(\gamma)$ & $e_0 $ & & & & & &\\
\hline
\end{tabular}
\caption{A basis $\mathcal{B}$ for $P(\mathbb{G}_2)$, and its dual $\mathcal{B}^*$.}
\end{table}

As in Section 7, these isomorphisms lift to the stable Auslander algebras of the representation-finite self-injective algebras of type $(\mathbb{D}_4,s,3)$ when $3 \nmid s$.  Then Theorem 4.2 of \cite{Per} allows us to deduce the following information about the periods of these algebras.
 
\begin{coro} Let $\Lambda$ be a representation-finite self-injective algebra of type $(\mathbb{D}_4, s, 3)$, and assume that $3 \nmid s$ and $s > 1$.  
\begin{enumerate}
\item If $\chr (k) = 2$, then $\Lambda$ has period $5s$.
\item If $\chr (k) \neq 2$, then $\Lambda$ has period $10s/(s,2)$. $\Box$
\end{enumerate}
\end{coro}

\section{Alternate approach: Triangulated orbit categories}
\setcounter{equation}{0}

We now describe another means of calculating the CY-dimensions and $\Sigma$-periods of $\Hom$-finite triangulated categories with finitely many indecomposables.  The basic idea, pursued originally by Holm and J\o rgensen \cite{HoJo1}, is to apply Amiot's characterization of such categories via covering theory \cite{Amiot}, while paying close attention to technical aspects of the covering theory of triangulated categories as demonstrated to us by Bernhard Keller \cite{KelP, TOCcor}.  In particular, this approach permits us to complement our previous calculations for the stable categories of representation-finite self-injective algebras, and we illustrate the process with several computations which eluded our previous methods.  In particular, we determine:
\begin{itemize}
\item The stable CY-dimensions of the representation-finite self-injective algebras of type $(\mathbb{A}_{2n+1},s,2)$ for even $n$.
\item The periods of the representation-finite self-injective algebras of type $(\mathbb{D}_n,s,2)$ when $n$ is odd and $s$ is even.
\item The periods of the representation-finite self-injective algebras of type $(\mathbb{E}_6,s,2)$ with $s$ even.
\end{itemize}
Similar arguments can also be used to compute the stable CY-dimensions and periods of most other representation-finite self-injective algebras, and even the CY-dimensions of algebraic triangulated categories with only finitely many indecomposables.  However, the current limitation of this approach occurs for triangulated categories with AR-quivers of the form $\mathbb{ZD}_{2n}/(\rho \tau^m)$ where $\rho$ has order $2$ (or $3$ when $n=2$).  The reason behind this difficulty is explained in the remarks preceding Propostion 9.7.

\vspace{3mm}
First of all, if $\mathcal{C}$ is a $k$-linear category and $F : \mathcal{C} \rightarrow \mathcal{C}$ is an automorphism, the orbit category $\mathcal{C}/F$ is defined to be the category with the same objects as $\mathcal{C}$, and with morphisms $$\mathcal{C}/F(x,y) = \bigoplus_{i \in \mathbb{Z}}\mathcal{C}(x,F^iy).$$  If $g : x \rightarrow F^i y$ and $f : y \rightarrow F^jz$, then the composite of $g$ and $f$ in $\mathcal{C}/F$ is given by $ F^i(f) \circ g : x \rightarrow F^{i+j}z$.  There is a natural projection functor $\pi : \mathcal{C} \rightarrow \mathcal{C}/F$, together with an isomorphism of functors $\eta : \pi \stackrel{\cong}{\longrightarrow} \pi F$ given by the identity maps $\eta_x = 1_{x}$ for each $x \in \mathcal{C}$.

Below, we want to consider orbit categories of a triangulated category $(\D, \Sigma)$ relative to an auto-equivalence $F$.   In order to define such an orbit category, we must first pass to a skeleton $\C$ of $\D$ and consider instead the orbit category of $\C$ with respect to the automorphism of $\C$ that corresponds to $F$ (Cf. \cite{Cover} Remark 3.8).  However, in this exchange, we must pay close attention to the triangulated structure of $\D$.  To make this treatment more precise, we let $I : \C \rightarrow \D$ be the inclusion functor and choose isomorphisms $\mu_x : x \rightarrow x_0$ with $x_0 \in \C$ for each $x \in D$ such that $\mu_x = 1_x$ for each $x \in \C$.  This defines a retraction $J : \D \rightarrow \C$ via $J(x) = x_0$ on objects and $J(f) = \mu_y f \mu_x^{-1}$ on morphisms $f \in \D(x,y)$, and we have $JI = 1_{\C}$ and $\mu : 1_{\D} \stackrel{\cong}{\longrightarrow} IJ$.  In order to make $\C$ into a triangulated category, we define a new suspension $\Sigma' := J\Sigma I : \C \rightarrow \C$, and define the distinguished triangles to be those isomorphic to a sequence of the form $$JA \stackrel{Jf}{\rightarrow} JB \stackrel{Jg}{\rightarrow} JC \stackrel{J(\Sigma(\mu_A) h)}{\longrightarrow} \Sigma'JA$$ where $A \stackrel{f}{\rightarrow} B \stackrel{g}{\rightarrow} C \stackrel{h}{\rightarrow} \Sigma A$ is a distinguished triangle in $\D$.  The axioms for triangulated categories (as stated in \cite{TCRTA}) can now be verified for $\C$ by using the fact that the functor $[I, \mu^{-1}_{\Sigma I} : I\Sigma' \rightarrow \Sigma I]$ takes a triangle in $\C$ to a triangle in $\D$, and the latter triangle is returned to the former by the functor $[J, J\Sigma(\mu) : J\Sigma \rightarrow \Sigma' J]$.  For instance, the rotation axiom (TR2) is checked via the following commutative diagram: 
$$\xymatrixcolsep{3.5pc} \xymatrix{JB \ar@{=}[d] \ar[r]^{Jg} & JC \ar@{=}[d] \ar[r]^{Jh} & J \Sigma A \ar[d]^{J\Sigma(\mu_A)} \ar[r]^{ -J\Sigma(\mu_b f)} & \Sigma'(JB) \ar@{=}[d] \\ JB \ar[r]_{Jg} & JC \ar[r]_{J(\Sigma(\mu_A) h)} & \Sigma'JA \ar[r]_{-\Sigma'Jf} & \Sigma'(JB)}$$
 The bottom row is the rotation of a typical triangle $JA \stackrel{Jf}{\rightarrow} JB \stackrel{Jg}{\rightarrow} JC \stackrel{J(\Sigma(\mu_A) h)}{\longrightarrow} \Sigma'JA$ in $\C$; while the top row is the image under $J$ of the rotation of the triangle $A \stackrel{f}{\rightarrow} B \stackrel{g}{\rightarrow} C \stackrel{h}{\rightarrow} \Sigma A$ in $\D$, and is hence a triangle in $\C$ by definition.  Since the vertical maps are all isomorphisms, the bottom row is also a triangle in $\C$.  With respect to this triangulated structure on $\C$ it is now immediate that the exact functors $[I,  \mu^{-1}_{\Sigma I}]$ and $[J, J\Sigma(\mu)]$ provide an equivalence of triangulated categories between $(\C, \Sigma')$ and $(\D, \Sigma)$.  If $F$ is an auto-equivalence of $\D$, then the corresponding automorphism of $\C$ is $F' = JFI$.  Moreover, if $S$ is a Serre functor for $\D$, then the automorphism $S' = JSI$ is easily checked to be a Serre functor for $\C$.

In general, if $\mathcal{C}$ is triangulated and $F$ is exact, $\mathcal{C}/F$ does not inherit a triangulated structure.  However, Keller demonstrated in \cite{TOC} that if $\mathcal{C}$ is triangle equivalent to the bounded derived category $D^b(\rmod k\Delta)$ with $\Delta$ an oriented Dynkin graph, and $F$ satisfies some mild hypotheses, then $\mathcal{C}/F$ is triangulated and $\pi : \mathcal{C} \rightarrow \mathcal{C}/F$ is exact.  In \cite{Amiot}, Amiot goes on to characterize the triangulated categories that arise in this way.   We have rephrased the following theorem in light of the above interpretation of orbit categories relative to auto-equivalences.   Moreover, from now on we will write $D^b(A)$ for the bounded derived category of $\rmod A$. 

\begin{therm}[Cf. \cite{Amiot}, Theorem 7.0.5] Suppose that $\mathcal{T}$ is a standard, algebraic Hom-finite triangulated category with only finitely many indecomposables up to isomorphism.  Then $\mathcal{T}$ is triangle equivalent to an orbit category $\D/F$ with $\D$ a skeleton of $D^b(\rmod kQ)$ for a Dynkin quiver $Q$ and an exact automorphism $F$ of $\D$.
\end{therm}

Clearly, the hypotheses of Amiot's theorem apply to $\mathcal{T} = \stmod \Lambda$ for a representation-finite standard self-injective $k$-algebra $\Lambda$.  Moreover, in order to find $A:=kQ$ and $F$, it suffices to match up the AR-quivers of $\mathcal{T}$ and $D^b(A)/F$--a task which can be easily accomplished using Happel's description of the AR-quivers of $D^b(A)$ and the effect of the suspension $\Sigma$ \cite{DCFDA}.  According to Happel, the AR-quiver of $D^b(k\Delta)$ is isomorphic to $\mathbb{Z}\Delta$.  Moreover, the category of indecomposables is realized by the mesh category: $\ind\ D^b(k\Delta) \approx k(\mathbb{Z}\Delta)$.  For the suspension functor $\Sigma$ on $D^b(k\Delta)$, we have $\Sigma^2 \cong \tau^{-h_{\Delta}}$.  Furthermore, $S = \Sigma \tau$ is a Serre functor on $D^b(k\Delta)$.  Finally the action of the suspension $\Sigma$ is given by the automorphism $\rho_{\Delta} \tau^{-\floor{h_{\Delta}/2}}$, where the automorphism $\rho_{\Delta}$ is the identity if $\Delta = \mathbb{D}_{2n}, \mathbb{E}_7$ or $\mathbb{E}_8$, and otherwise is as described in the definition of the $m$-fold mesh algebras in Section 3 (see Table I in \cite{DPic}).\\

\noindent
{\bf Examples.}  (1) Let $\mathcal{T} = \proj A$ where $A = P(\mathbb{L}_n)$ is the mesh algebra on $Q_{\mathbb{L}_n}$, and assume $\chr(k) \neq 2$.  The category $\mathcal{T} = \proj A$ is shown to be triangulated in \cite{Amiot}.  However, the suspension functor $\Sigma$ is not the identity as claimed there (due to what appears to be a misquoting of \cite{DPA}); rather, it is induced by the automorphism $\sigma$ of $A$ that fixes all vertices and all arrows other than $\epsilon$ in the quiver $Q_{\mathbb{L}_n}$, while sending $\epsilon \mapsto -\epsilon$ \cite{DPA}.  One easily checks that this functor is not isomorphic to the identity functor, or, equivalently, that the automorphism $\sigma$ of $A$ is not inner (except when $\chr (k) = 2$).  Indeed, the endomorphism ring of $Ae_0$ is generated as a $k$-algebra by $\epsilon$, and hence any endomorphism of $Ae_0$ commutes with $\epsilon$, implying that there is no endomorphism $\eta_{Ae_0}$ such that $-\epsilon \eta_{Ae_0} \cong \eta_{Ae_0} \epsilon$.  Since $\mathcal{T}$ is standard and algebraic, it can be realized as an orbit category of some $D^b(kQ)$ with $Q$ Dynkin.  Since the AR-quiver of $\mathcal{T}$ coincides with the bound quiver of $P(\mathbb{L}_n)$, as a translation quiver it is $\mathbb{ZA}_{2n}/\gen{\rho \tau}$ where $\rho =\rho_{\mathbb{A}_{2n}}$ has order $2$.  Hence, as a triangulated category, $\mathcal{T}$ is equivalent to $\D/F$ where $\D$ is a skeleton of $D^b(k\mathbb{A}_{2n})$ and $F = \tau_{\D}^n \Sigma_{\D}$ (The AR-quiver of $\D$ is $\mathbb{ZA}_{2n}$ and the automorphism induced by $\Sigma_{\D}$ translates $n + \frac{1}{2}$ units to the right and reflects in the center line: see eg. \cite{DPic}).

In 7.4 of \cite{TOC}, Keller claims\footnote{This claim has since been corrected \cite{TOCcor}.} that the isomorphism of functors $\tau_{\D}^{-(2n+1)} \cong \Sigma_{\D}^2$ in $\D$ induces the following isomorphisms in the orbit category $\mathcal{T}$:
\begin{equation} 1_{\mathcal{T}} \cong F^2 = (\tau_{\mathcal{T}}^n \Sigma_{\mathcal{T}})^2 \cong \tau_{\mathcal{T}}^{2n} \Sigma_{\mathcal{T}}^2 \cong \tau_{\mathcal{T}}^{-1}.\end{equation}
However, this contradicts what we observed above: $\tau_{\T}^{-1} \cong \tau_{\T} \cong \Sigma_{\T} \ncong 1_{\mathcal{T}}$.  See the Remark following Proposition 9.6 for a discussion of the error in (9.1).  Moreover, on one hand, we would like to say the functor induced by $F$ is $\tau_{\T}^n \Sigma_{\T}$, while on the other hand, the functor induced by $F$ in $\D/F$ should be the identity.  Nevertheless, observe that $\tau_{\T}^n \Sigma_{\T} \cong \tau_{\T}^{n+1} \ncong 1_{\mathcal{T}}$ whenever $n$ is even.\\

\noindent
(2) For another example, let $\mathcal{T} = \stmod \Lambda$ where $\Lambda$ is the representation-finite self-injective algebra of type $(\mathbb{D}_6, 5/3,1)$.  We assume that $\chr(k) \neq 2$.  Since the stable AR-quiver of $\Lambda$ is isomorphic to $\mathbb{Z}\mathbb{D}_6/\gen{\tau^{15}}$, $\mathcal{T}$ is triangle equivalent to the orbit category $\D/F$ where $\D$ is a skeleton of $D^b(k \mathbb{D}_6)$ and $F \cong \tau^{15}$ (we fix any orientation for $\mathbb{D}_6$).  The AR-quiver of $\D$ is isomorphic to $\mathbb{Z}\mathbb{D}_6$, and the suspension $\Sigma_{\D}$ acts on the AR-quiver by shifting everything $5$ units to the right: $\Sigma_{\D} \cong \tau_{\D}^{-5}$.  The Serre functor on $\D$ is $S_{\D} = \Sigma_{\D} \tau_{\D} \cong \tau_{\D}^{-4}$.  In particular, we see that we can take $F = \Sigma_{\D}^3$, which yields functorial isomorphisms
 \begin{equation} \Sigma_{\T}^3 \pi \cong \pi \Sigma_{\D}^3 = \pi F \cong \pi, \end{equation}
as the covering functor $\pi : \D \rightarrow \D/F$ is exact with $\pi \cong \pi F$.  Since $\pi$ is dense, it follows that $\Sigma_{\T}^3$ and $1_{\T}$ have the same effect on the objects of $\T$.  However, this isomorphism does not guarantee that $\Sigma_{\T}^3 \cong 1_{\T}$ as $\pi$ is not full.  In fact, since $\Lambda$ is Schurian, an isomorphism $\Sigma_{\T}^3 \cong 1_{\T}$ would imply that $\Omega^3_{\Lambda^e}(\Lambda) \cong \Lambda$ (by Lemma 4.4 of \cite{Per}), which contradicts the calculation in \cite{Per} that the (minimal) period of $\Lambda$ is $6$.\\

To avoid the potential pitfalls illustrated above, we need to consider some extra structure on endofunctors of $D^b(A)$.  Namely, in order to obtain well-defined induced endofunctors on the orbit category in which we are interested, we need to consider $F$-equivariant functors on $D^b(A)$.  In particular, we need to keep track of {\it how} these endofunctors commute with $F$.  To provide more detail, we now review some basic elements of the theory recently introduced by Asashiba in \cite{Cover}.

\begin{defin}[Cf. \cite{Cover}]  Let $\C$ be a $k$-linear category with an automorphism $F: \C \rightarrow \C$.  A {\bf (right) $F$-invariant functor} from $\C$ to $\C'$ is a pair $(H,\phi)$ where $H : \C \rightarrow \C'$ is a $k$-linear functor and $\phi : H \rightarrow HF$ is an isomorphism.  A {\bf (weakly) $F$-equivariant functor} on $\C$ is a pair $(E, \psi)$ where $E : \C \rightarrow \C$ is a $k$-linear functor and $\psi : FE \rightarrow EF$ is an isomorphism.
\end{defin}

We point out that these terms are defined more generally in \cite{Cover}, with respect to a group $G$ of automorphisms of $\C$.  For our purposes, $G = \gen{F}$ is always cyclic and it thus suffices to specify only one functorial isomorphism in each case.  A morphism $\alpha : (H,\phi) \rightarrow (H',\phi')$ of $F$-invariant functors  is a morphism $\alpha : H \rightarrow H'$ satisfying $\alpha_{Fx} \phi_x = \phi'_x \alpha_x$ for all $x \in \C$.  Similarly, a morphism $\alpha : (E, \psi) \rightarrow (E',\psi')$ of $F$-equivariant functors is a morphism $\alpha : E \rightarrow E'$ satisfying $\psi'_x \circ F(\alpha_x) = \alpha_{Fx} \circ \psi_x$ for all $x \in \C$.  We shall denote the categories of $F$-invariant functors from $\C$ to $\C'$ and of $F$-equivariant functors on $\C$ by $\mathrm{inv}_F(\C,\C')$ and $\mathrm{equ}_F(\C,\C)$ respectively.

Observe that the natural covering functor $\pi : \C \rightarrow \C/F$ together with the isomorphism $\eta$ described above is an $F$-invariant functor.  Furthermore, the orbit category $\C/F$ is characterized by the following universal property (cf. 2.6, 2.7 in \cite{Cover}): \begin{quote} \em For every $F$-invariant functor $(H,\phi) : \C \rightarrow \C'$ there exists a unique functor $G : \C/F \rightarrow \C'$ (up to isomorphism) such that $(H,\phi) \cong (G \pi, G \eta)$ as $F$-invariant functors.  Equivalently, composition with $(\pi, \eta)$ induces an isomorphism of categories from the category $\mathrm{fun}_k(\C/F,\C')$ of $k$-linear functors $\C/F \rightarrow \C'$ to the category $\mathrm{inv}_F(\C,\C')$ of $F$-invariant functors $\C \rightarrow \C'$.\end{quote}  

Concerning composition, notice that the composite of two $F$-equivariant functors $(E,\psi)$ and $(E',\psi')$ yields another $F$-equivariant functor $(E E', E(\psi') \circ \psi_{E'})$.  Furthermore, if $(H,\phi) : \C \rightarrow \C'$ is an $F$-invariant functor, then $HE$ will inherit a natural $F$-invariant structure as well.  To be precise, $(HE, H(\psi)\phi_E) : \C \rightarrow \C'$ is $F$-invariant, and moreover one can check that these compositions are associative.  In particular, taking $(H,\phi) = (\pi, \eta)$, we get functors $$\mathrm{equ}_F(\C,\C) \stackrel{(\pi,\eta)\circ -}{\longrightarrow} \mathrm{inv}_F(\C,\C/F) \stackrel{\cong}{\longrightarrow} \mathrm{fun}_k(\C/F,\C/F)$$ which respect composition.  Via this sequence of maps, we see that an $F$-equivariant functor $(E,\psi) : \C \rightarrow \C$ {\it induces} a functor $\bar{E} : \C/F \rightarrow \C/F$.  More precisely, the action of $\bar{E}$ on objects of $\C/F$ coincides with that of $E$ on $\C$, while $\bar{E}$ takes a map $f : x \rightarrow Fy$ (regarded as a map in $\C/F$) to the composite $$Ex \stackrel{Ef}{\longrightarrow} EFy \stackrel{\psi_y^{-1}}{\longrightarrow} FEy.$$

Of special importance are the $F$-equivariant functors of the form $(1_{\C},\epsilon)$ where $\epsilon: F \rightarrow F$ is an isomorphism.  We shall denote the induced functor $\Delta(\epsilon) : \C/F \rightarrow \C/F$.  To illustrate, consider the right $F$-equivariant functors $(F, 1_{F^2})$ and $(F,-1_{F^2})$.  The composite $(\pi,\eta)(F,1_{F^2}) = (\pi F, \eta_F)$ is isomorphic to $(\pi, \eta)$ as an $F$-invariant functor via $\eta : \pi \rightarrow \pi F$.  It follows that $(F, 1_{F^2})$ induces the identity functor $1_{\C/F}$.  On the other hand $(F, -1_{F^2})$ is typically not isomorphic to $(F,1_{F^2})$ as an $F$-equivariant functor, and hence will not necessarily induce the identity on $\C/F$.  In fact, as $(F, -1_{F^2}) = (1_{\C},-1_F) \circ (F,1_{F^2})$, we see that the functor induced by $(F,-1_{F^2})$ will coincide with $\Delta(-1_F)$.

\begin{lemma}[Keller \cite{TOC, TOCcor}]  Let $[F,\alpha] : \D \rightarrow \D$ be an exact automorphism of the triangulated category $(\D,\Sigma)$ as in \cite{TOC}, where $\alpha : F\Sigma \stackrel{\cong}{\longrightarrow} \Sigma F$.  Then $(\Sigma, \alpha)$ is $F$-equivariant and induces the suspension $\Sigma_{\D/F}$ of the orbit category $\D/F$. $\Box$
\end{lemma}

Notice that if we take $[F,\alpha] = [\Sigma, -1_{\Sigma^2}]$, then the above remarks show that the suspension of the orbit category $\D/\Sigma$ is isomorphic to $\Delta(-1_{\Sigma})$, which is not isomorphic to the identity in general.  Likewise, if we factor out $F = \Sigma^3$ as in Example 2, we see that the functor induced on the orbit category by $(\Sigma^3, -1_{\Sigma^4})$ is $\Sigma_{\D/F}^3 \cong \Delta(-1_{\D/F})$ rather than the identity.

\vspace{3mm}
We now consider the relationship between the Serre functors of $\C$ and $\C/F$.  We assume that both $\C$ and $\C/F$ have finite-dimensional Hom-spaces, which is certainly true when $\C$ and $F$ are as in the main theorem of \cite{TOC}.  Recall from \cite{CYTC} that a Serre functor $S : \C \rightarrow \C$ is characterized by {\it trace maps}, which are functorial isomorphisms $$t_x : \C(-,Sx) \rightarrow D \C(x,-)$$ natural in $x \in \C$.  For any automorphism $F$ of $\C$, we can define a unique isomorphism $\sigma_F : SF \rightarrow FS$ so that the following diagram commutes for any $x \in \C$.

\begin{equation} \vcenter{\xymatrixcolsep{3.0pc} \xymatrix{ \C(-,SFx) \ar[r]^{\C(-,\sigma_{F,x})} \ar[d]^{t_{Fx}} & \C(-,FSx) \ar[r]^{F^{-1}} & \C(F^{-1}(-), Sx) \ar[d]^{(t_x)_{F^{-1}}} \\ D \C(Fx,-) \ar[rr]^{DF} & & D\C(x,F^{-1}(-))}} \end{equation}

We recall that there is a unique way (up to isomorphism) to make a Serre functor $S$ into a triangulated functor $[S,\eta]$ such that the diagram
 
 \begin{equation} \vcenter{\xymatrixcolsep{3.0pc} \xymatrix{ \C(\Sigma(-),S\Sigma x) \ar[d]^{(t_{\Sigma x})_{\Sigma}} \ar[r]^{\C(\Sigma(-),\eta_x)} & \C(\Sigma(-),\Sigma Sx) \ar[r]^{\Sigma^{-1}} & \C(-, Sx) \ar[d]^{t_x} \\ D \C(\Sigma x, \Sigma(-)) \ar[rr]^{-D\Sigma} & & D\C(x,-)}} \end{equation}
 commutes for all $x \in \C$ \cite{GCY3, CYTC, CYTCcor}.

\begin{propos} Suppose that $\C$ and $\C/F$ are Hom-finite triangulated categories and $S$ is a Serre functor on $\C$.
\begin{enumerate}
\item[(a)] The $F$-equivariant functor $(S,\sigma_F^{-1})$ induces the Serre functor $\overline{S}$ of $\C/F$.
\item[(b)] We have $\sigma_S = 1_{S^2}$.
\item[(c)] The (unique) enhancement of $S$ into a triangle functor is given by $[S,-\sigma_{\Sigma}]$.
\item[(d)] If $F = F_1 F_2$ for $F_1, F_2 \in \aut (\C)$, then $\sigma_F = F_1(\sigma_{F_2}) \circ (\sigma_{F_1})_{F_2}$.
\end{enumerate}
\end{propos}

\noindent
{\it Proof.}  (a)  For $x \in \C$ we define trace maps $$t_{\pi x} : \C/F(-,\overline{S}(\pi x)) = \oplus_{p \in \mathbb{Z}} \C(-,F^pSx) \rightarrow \oplus_{p \in \mathbb{Z}} D\C(x,F^{-p}(-)) = D(\C/F(\pi(x),-))$$ via the composites $$\C(-,F^pSx) \stackrel{\C(-,(\sigma^{-1}_{F^p})_x)}{\longrightarrow} \C(-,SF^p x) \stackrel{t_{F^px}}{\longrightarrow} D\C(F^p x,-) \stackrel{D(F^p)}{\longrightarrow} D\C(x,F^{-p}(-)).$$  One can now use the definition of the induced functor $\overline{S}$ to check that these trace maps are natural in $x$. 

(b) For any $x,y \in \C$, two applications of Lemma I.1.1 of \cite{Noeth} yield the two commutative triangles below.
$$\xymatrix{\C(y,Sx) \ar[r]^{t_{x,y}} \ar[d]^S \ar[dr]^{Dt_{y,Sx}} & D\C(x,y) \\ \C(Sy,S^2x) \ar[r]_{t_{Sx,Sy}} & D(Sx,Sy) \ar[u]_{DS}}$$
Comparing the square above to the diagram used to define $\sigma_S$ now reveals that $\sigma_S = 1_{S^2}$.

(c) If the triangulated enhancement of $S$ is $[S,\eta]$, the difference in sign between (9.3) and (9.4) implies that $\eta = -\sigma_{\Sigma}$.

(d) This can be seen by pasting together the natural transformations in the diagram
$$\vcenter{\xymatrix{\C \ar[r]^{F_2} \ar[d]_S & \C \ar[r]^{F_1} \ar[d]_S \ar@{=>}[dl]|{\sigma_{F_2}} & \C \ar[d]^S \ar@{=>}[dl]|{\sigma_{F_1}} \\ \C \ar[r]_{F_2} &\C \ar[r]_{F_1} & \C}}.\ 
\Box$$ \\

For $p \in \mathbb{N}$, part (d) shows that $\sigma_{F^p}$ is the isomorphism $SF^p \rightarrow F^pS$ obtained by pasting together the natural transformations in the following diagram, in which $F$ appears $p$ times in each row.

$$\xymatrix{\C \ar[r]^{F} \ar[d]_S & \C \ar[r]^{F} \ar[d]_S \ar@{=>}[dl]|{\sigma_{F}} & \cdots \ar[r]^F \ar@{=>}[dl]|{\sigma_{F}} & \C \ar[d]^S \ar@{=>}[dl]|{\sigma_{F}} \\ \C \ar[r]_{F} &\C \ar[r]_{F} & \cdots \ar[r]_F & \C}$$
We shall also write $\sigma_{F}^{(m)}$ for the isomorphism $S^m F \rightarrow F S^m$ obtained by pasting together the natural transformations in the diagram below, in which $S$ appears $m$ times in each row.
$$\xymatrix{ \C \ar[d]_{F} \ar[r]^S & \C \ar[d]_{F} \ar[r]^S & \cdots \ar[r]^S & \C \ar[d]^{F} \\ \C \ar[r]_S \ar@{=>}[ur]|{\sigma_{F}} & \C \ar[r]_S  \ar@{=>}[ur]|{\sigma_{F}} & \cdots \ar[r]_S  \ar@{=>}[ur]|{\sigma_{F}} & \C}
$$

The following theorem and proof are essentially due to Keller, who supplied the argument when $m=1$ \cite{KelP, TOCcor}.
 
\begin{therm}  Let $\C$ be a skeleton of $D^b(A)$ with suspension $\Sigma$ and Serre functor $S$, which we may assume are automorphisms, and set $F = S^m \Sigma^d$ for $m, d \in \mathbb{Z}$.  Then $(\Sigma, (-1)^{d+m}(\sigma_{\Sigma}^{(m)})_{\Sigma^d})$ and $(S, S^m(\sigma_{\Sigma^d}^{-1}))$ are $F$-equivariant functors inducing the suspension $\Sigma_{\C/F}$ and Serre functor $S_{\C/F}$, respectively, on the orbit category $\C/F$.  In particular, we have $$S_{\C/F}^m \Sigma_{\C/F}^d \cong \Delta((-1)^{d(d+m)}).$$
\end{therm}

\noindent
{\it Proof.}  We need to identify the correct $F$-equivariant structures on $S$ and $\Sigma$.  By Lemma 9.3, the $F$-equivariant structure on $\Sigma$ comes from the isomorphism $\alpha : F\Sigma \rightarrow \Sigma F$ that makes $F$ into a triangle functor.  As a triangle functor $F = [S,-\sigma_{\Sigma}]^m [\Sigma^d, (-1)^d \cdot 1_{\Sigma^{d+1}}] \cong [S^m \Sigma^d, (-1)^{d+m}(\sigma_{\Sigma}^{(m)})_{\Sigma^d}]$.  In fact, the isomorphism $\alpha$ here can be seen by pasting together the left-most column of squares in the figure below.  Thus $\Sigma_{\C/F}$ is induced by $(\Sigma, (-1)^{d+m}(\sigma_{\Sigma}^{(m)})_{\Sigma^d})$.  For $S$ we need to find $\sigma_F : S^{m+1}\Sigma^d \rightarrow S^m \Sigma^d S$.  By Proposition 9.4(b), $\sigma_{S^m} = 1_{S^{m+1}}$ and hence $\sigma_F = S^m(\sigma_{\Sigma^d})$ by Proposition 9.4(d).  In the figure below, $\sigma_F^{-1}$ is obtained by pasting the squares in the central column.

Now we consider the $F$-equivariant structure on $F$ obtained by writing $F = (S,\sigma_F^{-1})^m \circ (\Sigma, \alpha)^d$.  The corresponding automorphism $\epsilon$ of $F^2$ is then given as the composite $$FS^m\Sigma^d \stackrel{(\sigma_F^{(-m)})_{\Sigma^d}}{\longrightarrow} S^mF\Sigma^d \stackrel{S^m(\alpha^{(d)})}{\longrightarrow} S^m\Sigma^d F,$$ where $\sigma_F^{(-m)}$ and $\alpha^{(d)}$ are shorthand for the isomorphisms $FS^m \cong S^m F$ and $F\Sigma^d \cong \Sigma^d F$ induced by $\sigma_F^{-1}$ and $\alpha$ respectively.  This construction of $\epsilon$ is also illustrated in the diagram below, in which $\alpha^{(d)}$ is obtained by pasting the squares on the left half and $\sigma_F^{(-m)}$ is obtained by pasting the squares on the right half.  In each row and column of the diagram, $\Sigma$ appears $d$ times and $S$ appears $m$ times.  Note also that we have abbreviated $\sigma_{\Sigma}$ by $\sigma$.

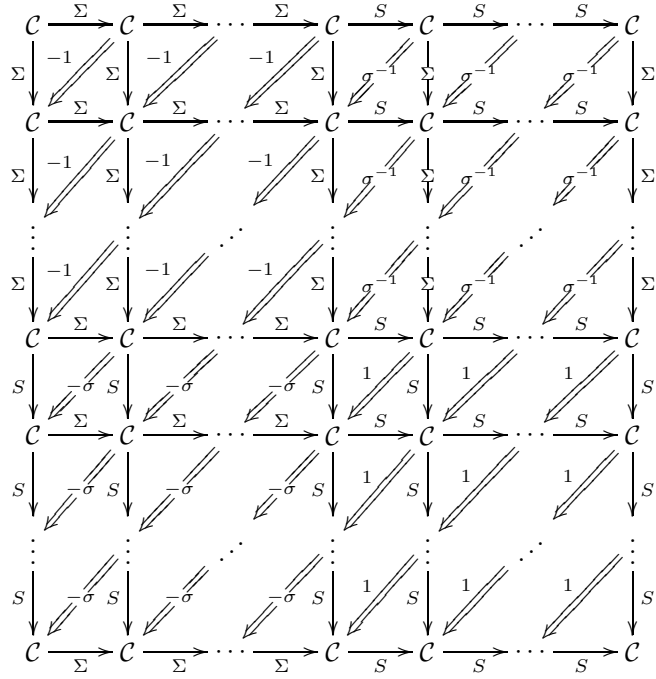
$$\xymatrix{ \C\ar[r]^{\Sigma} \ar[d]_{\Sigma} & \C \ar[r]^{\Sigma} \ar[d]_{\Sigma} \ar@{=>}[dl]_{-1} & \cdots  \ar[r]^{\Sigma}  \ar@{=>}[dl]_{-1} & \C  \ar[r]^{S} \ar[d]_{\Sigma} \ar@{=>}[dl]_{-1} & \C  \ar[r]^{S} \ar[d]|{\Sigma} \ar@{=>}[dl]|{\sigma^{-1}} & \cdots  \ar[r]^{S} \ar@{=>}[dl]|{\sigma^{-1}} & \C \ar[d]^{\Sigma} \ar@{=>}[dl]|{\sigma^{-1}} \\ 
 \C\ar[r]^{\Sigma} \ar[d]_{\Sigma} & \C \ar[r]^{\Sigma} \ar[d]_{\Sigma} \ar@{=>}[dl]_{-1} & \cdots  \ar[r]^{\Sigma}  \ar@{=>}[dl]_{-1} & \C  \ar[r]^{S} \ar[d]_{\Sigma} \ar@{=>}[dl]_{-1} & \C  \ar[r]^{S} \ar[d]|{\Sigma} \ar@{=>}[dl]|{\sigma^{-1}} & \cdots  \ar[r]^{S} \ar@{=>}[dl]|{\sigma^{-1}} & \C \ar[d]^{\Sigma} \ar@{=>}[dl]|{\sigma^{-1}} \\ 
\vdots \ar[d]_{\Sigma} & \vdots \ar[d]_{\Sigma} \ar@{=>}[dl]_{-1} & \adots \ar@{=>}[dl]_{-1} & \vdots \ar[d]_{\Sigma} \ar@{=>}[dl]_{-1} & \vdots \ar[d]|{\Sigma} \ar@{=>}[dl]|{\sigma^{-1}} & \adots  \ar@{=>}[dl]|{\sigma^{-1}}  & \vdots \ar@{=>}[dl]|{\sigma^{-1}} \ar[d]^{\Sigma} \\
\C\ar[r]^{\Sigma} \ar[d]_{S} & \C \ar[r]^{\Sigma} \ar[d]_{S} \ar@{=>}[dl]|{-\sigma} & \cdots  \ar[r]^{\Sigma}  \ar@{=>}[dl]|{-\sigma} & \C  \ar[r]^{S} \ar[d]_{S} \ar@{=>}[dl]|{-\sigma} & \C  \ar[r]^{S} \ar[d]_{S} \ar@{=>}[dl]_{1} & \cdots  \ar[r]^{S} \ar@{=>}[dl]_{1} & \C \ar[d]^{S} \ar@{=>}[dl]_{1} \\ 
\C\ar[r]^{\Sigma} \ar[d]_{S} & \C \ar[r]^{\Sigma} \ar[d]_{S} \ar@{=>}[dl]|{-\sigma} & \cdots  \ar[r]^{\Sigma}  \ar@{=>}[dl]|{-\sigma} & \C  \ar[r]^{S} \ar[d]_{S} \ar@{=>}[dl]|{-\sigma} & \C  \ar[r]^{S} \ar[d]_{S} \ar@{=>}[dl]_{1} & \cdots  \ar[r]^{S} \ar@{=>}[dl]_{1} & \C \ar[d]^{S} \ar@{=>}[dl]_{1} \\ 
\vdots \ar[d]_{S} & \vdots \ar[d]_{S} \ar@{=>}[dl]|{-\sigma} & \adots \ar@{=>}[dl]|{-\sigma} & \vdots \ar[d]_{S} \ar@{=>}[dl]|{-\sigma} & \vdots \ar[d]_{S} \ar@{=>}[dl]_{1} & \adots  \ar@{=>}[dl]_{1}  & \vdots \ar@{=>}[dl]_{1} \ar[d]^{S} \\
\C \ar[r]_{\Sigma} & \C \ar[r]_{\Sigma} & \cdots \ar[r]_{\Sigma} & \C \ar[r]_{S} & \C \ar[r]_{S} & \cdots \ar[r]_{S} & \C}$$

In particular, notice that $\alpha^{(d)}$ coincides with $(-1)^{d(d+m)}(\sigma^{(m)}_{\Sigma^d})_{\Sigma^d}$.  It follows that $$\epsilon = S^m((-1)^{(d+m)d}\sigma_{\Sigma^d}^{(m)})_{\Sigma^d} \circ S^m(\sigma_{\Sigma^d}^{(-m)})_{\Sigma^d} = (-1)^{(d+m)d} \cdot 1_{F^2}.$$  Consequently the functor induced on $\C/F$ by $(F,\epsilon) = (F, 1_{F^2}) \circ (1, \epsilon)$ is isomorphic to both $\Delta(\epsilon) = \Delta((-1)^{d(d+m)})$ and to $S_{\C/F}^m \Sigma_{\C/F}^d$.  $\Box$ \\

The above theorem is critical for determining when relations between the automorphisms $S$ and $\Sigma$ on $\C$ continue to hold in the orbit category. 
We now illustrate how it may be applied to determine the CY-dimensions and $\Sigma$-periods of the finite type stable categories $\stmod \Lambda$.  We begin by verifying the stable Calabi-Yau dimensions computed for the Moebius algebras in \cite{HoJo1}.

\begin{propos}[Cf. 7.1; \cite{HoJo1}, Theorem 5.2] Let $\Lambda$ be a representation-finite self-injective algebra of type $(\mathbb{A}_{2n+1},s,2)$.  Then $\stmod \Lambda$ is Calabi-Yau if and only if $(n+1,s)=1$ and $n \equiv s \ (\mbox{mod}\ 2)$.  In this case the stable CY-dimension of $\Lambda$ is $d = K_{n,s}(2n+1)-1$ where $$K_{n,s} := \inf\ \left\{ r\ \left| \ r \geq 1, r(n+1) \equiv 1\ (\mbox{mod}\ s),\ 2|\frac{r(s+n+1)-1}{s} \right. \right\}.$$ 
\end{propos}

\noindent
{\it Proof.}  As the necessary conditions were already encountered in 7.1, we need only verify the claimed CY-dimension.  We thus assume that $n \equiv s\ (\mbox{mod}\ 2)$ and $(n+1,s)=1$.  Recall that the AR-quiver of $\Lambda$ is isomorphic to $\mathbb{ZA}_{2n+1}/\gen{\rho \tau^{s(2n+1)}}$ with $\rho^2 =1$.  Moreover, on a skeleton $\D$ of $D^b(\rmod k\mathbb{A}_{2n+1})$ we have isomorphisms $\Sigma_{\D} \cong \rho \tau^{-n-1}$ and $S_{\D} \cong \Sigma_{\D} \tau \cong \rho \tau^{-n}$.  Thus $\stmod \Lambda$ is realized as the orbit category $\D/F$ where $F = \Sigma_{\D}^{-s-n}S_{\D}^{n+1-s}$.  Since $n$ and $s$ have the same parity, Theorem 9.5 implies that $\Sigma^{s+n} \cong S^{n+1-s}$ for the suspension $\Sigma$ and Serre functor $S$ on $\stmod \Lambda$.  Likewise, since we have an isomorphism $\Sigma^{2n}_{\D} \cong S^{2(n+1)}_{\D}$, we can realize $\stmod \Lambda$ as the orbit category $\D/F'$ with $F' = \Sigma_{\D}^{n-s}S_{\D}^{-s-n-1}$.  Again by Theorem 9.5 we conclude that $\Sigma^{n-s} \cong S^{n+s+1}$, and combining these isomorphisms yields $\Sigma^{2n} \cong S^{2(n+1)}$ on $\stmod \Lambda$.  Now notice that by our choice of $K_{n,s}$, we get integers $$a := \frac{K_{n,s}(n+1)-1}{s}\ \mbox{and}\ b := \frac{K_{n,s}(s-n-1)+1}{2s}$$ satisfying $(n+1-s)a + 2(n+1)b = 1$.  Consequently $$S = (S^{n+1-s})^a(S^{2(n+1)})^b \cong (\Sigma^{s+n})^a(\Sigma^{2n})^b = \Sigma^{K_{n,s}(2n+1)-1}. \Box $$ 

\vspace{2mm}
\noindent
{\bf Remark.}  Returning briefly to Example (1) we can now identify the error in (9.1). In fact, it lies in the last isomorphism $\tau_{\T}^{2n}\Sigma_{\T}^2 \cong \tau_{\T}^{-1}$ as the isomorphism $\tau^{-(2n+1)}_{\D} \cong \Sigma_{\D}^2$ fails to pass down to the orbit category.  Indeed, using the identity $S = \Sigma \tau$ we convert this isomorphism to $S_{\D}^{2n+1} \cong \Sigma_{\D}^{2n-1}$, and attempt to mimic the argument in the above proof by applying Theorem 9.5 for two different choices of $F$.  If $F_1 = \tau_{\D}^n \Sigma_{\D} \cong S_{\D}^n\Sigma_{\D}^{1-n}$ and $F_2 = S_{\D}^{-(1+n)} \Sigma_{\D}^{n}$, then the functors induced by $F_1$ and $F_2$ will differ by $\Delta(-1)$, yielding $\tau_{\T}^{-(2n+1)} \cong \Delta(-1) \Sigma_{\T}^2$.  In particular, $\tau_{\T} \cong \Sigma_{\T} \cong \Delta(-1)$ here.\\

While we would like to similarly apply Theorem 9.5 to compute the CY-dimensions for the stable categories of algebras of type $(\mathbb{D}_n, s, 2)$ when $s$ is even, in this case the necessary automorphism $F$ of $\D \approx D^b(\rmod k\mathbb{D}_n)$ cannot be expressed in terms of $S_{\D}$ and $\Sigma_{\D}$ alone.  Indeed, we would need $F = \rho \tau^{s(2n-3)}$, but the necessary condition $(n-1,s) = 1$ forces $n$ to be even, making $\Sigma_{\D} \cong \tau_{\D}^{1-n}$ and $S_{\D} \cong \tau_{\D}^{2-n}$.  Thus we encounter a need to describe precisely how the degree-$2$ automorphism $\rho$ of $\D$ commutes with $\Sigma_{\D}$ and $S_{\D}$ -- a task which we do not pursue here.  The same problem occurs for the triangulated categories with AR-quivers of the form $\mathbb{ZD}_4/\gen{\rho\tau^m}$ where $\rho$ is a degree-$3$ automorphism of $\mathbb{ZD}_4$.  We can however use Theorem 9.5 to obtain the periods of algebras of type $(\mathbb{D}_n, s, 2)$ when $s$ is even and $n$ is odd (Cf. Proposition 7.3).

\begin{propos} Let $\Lambda$ be a representation-finite self-injective algebra of type $(\mathbb{D}_{n},s,2)$ with $n$ odd and $s \geq 2$ even.  Then the period of $\Lambda$ is 
$$ p_{\Lambda} = \left\{ \begin{array}{rl} \frac{s(2n-3)}{(s,n-1)}, & \mbox{if\ } \chr(k) = 2\ \mbox{and}\  2\ |\ \frac{s+n-1}{(s,n-1)}\\ \frac{2s(2n-3)}{(s,n-1)}, & \mbox{otherwise.}
 \end{array} \right.$$
\end{propos}

\noindent
{\it Proof.}  Since the AR-quiver of $\Lambda$ has shape $\mathbb{ZD}_n/\gen{\rho \tau^{s(2n-3)}}$, we can realize $\stmod \Lambda$ as the orbit category $\D/F$ where $\D$ is a skeleton of $D^b(\rmod k\mathbb{D}_n)$ and $F = \rho \tau_{\D}^{s(2n-3)}$.  The identities $\Sigma_{\D} \cong \rho \tau_{\D}^{1-n}$ and $S_{\D} \cong \rho \tau_{\D}^{2-n}$ (for $n$ odd) imply that we can take $F = S^{n-s-1}_{\D} \Sigma_{\D}^{-s-n+2}$.  By Theorem 9.5, we now have $\Sigma^{s+n-2} \cong \Delta(-1) S^{n-s-1}$ for the Serre functor $S$ and suspension $\Sigma$ of $\stmod \Lambda$.  As in the previous proof, the identity $\Sigma_{\D}^{2(2-n)} \cong S_{\D}^{2(1-n)}$ can be carried down to an isomorphism $\Sigma^{2(2-n)} \cong S^{2(1-n)}$ on the orbit category.  Setting $$a : = \frac{n-s-1}{(n-s-1,2(n-1))}\ \mbox{and}\ b := \frac{2(1-n)}{(n-s-1,2(n-1))}$$ (with the convention that the g.c.d $(x,y)$ of two integers is always positive), we obtain $$\Sigma^{2(n-2)a + (s+n-2)b} = (\Sigma^{2(2-n)})^{-a} (\Sigma^{s+n-2})^b \cong (S^{2(1-n)})^{-a} (\Delta(-1)S^{n-s-1})^b = \Delta(-1)^b.$$  Simplifying the exponent of $\Sigma$ we obtain $$\Sigma^{\frac{2s(2n-3)}{(n-1-s, 2(n-1))}} \cong \Delta(-1)^{2(n-1)/(n-s-1, 2(n-1))} = \left\{ \begin{array}{ll} 1_{\stmod \Lambda}, & \mbox{if\ } 2 \nmid \frac{s+n-1}{(s,n-1)} \\ \Delta(-1), & \mbox{if\ } 2 \mid \frac{s+n-1}{(s,n-1)}. 
\end{array} \right.$$
Notice that in the first case we have $(n-1-s, 2(n-1)) = (n-1-s, n-1) = (s,n-1)$, while in the second case $(n-1-s,2(n-1)) = 2(s,n-1)$.  Since $\Delta(-1) \cong 1_{\stmod \Lambda}$ if and only if $\chr(k)=2$, the stated periods for $\Lambda$ now follow. $\Box$\\

Finally, we determine the periods of the remaining representation-finite self-injective algebras of type $(\mathbb{E}_6,s,2)$ (Cf. Proposition 7.4).

\begin{propos} Let $\Lambda$ be a representation-finite self-injective algebra of type $(\mathbb{E}_6,s,2)$ with $s \geq 2$ even.  Then the period of $\Lambda$ is 
$$ p_{\Lambda} = \left\{ \begin{array}{rl} \frac{11s}{(s,6)}, & \mbox{if\ } \chr(k) = 2\ \mbox{and}\  s \equiv 2\ (\mbox{mod}\ 4) \\ \frac{22s}{(s,6)}, & \mbox{otherwise.}
 \end{array} \right.$$
\end{propos}

\noindent
{\it Proof.}  The proof proceeds completely parallel to the previous proof so we only outline the details.  We can realize $\stmod \Lambda$ as the orbit category $\D/F$ where $\D$ is a skeleton of $D^b(\rmod k\mathbb{E}_6)$ and $F = \rho \tau_{\D}^{-11s}$.  The identities $\Sigma_{\D} \cong \rho \tau_{\D}^{-6}$ and $S_{\D} \cong \rho \tau_{\D}^{-5}$  imply that we may actually select $F = S^{s-6}_{\D} \Sigma_{\D}^{s+5}$.  By Theorem 9.5, we now have $\Sigma^{s+6} \cong \Delta(-1) S^{6-s}$ for the Serre functor $S$ and suspension $\Sigma$ of $\stmod \Lambda$.  Once again, the identity $\Sigma_{\D}^{10} \cong S_{\D}^{12}$ can be carried down to an isomorphism $\Sigma^{10} \cong S^{12}$ on $\stmod \Lambda$.  Using these identities we find that $$\Sigma^{22s/(12,s-6)} \cong \Delta(-1)^{12/(12,s-6)} = \left\{ \begin{array}{ll} 1_{\stmod \Lambda}, & \mbox{if\ } s \equiv 0 \ (\mbox{mod}\ 4) \\ \Delta(-1), & \mbox{if\ } s \equiv 2\ (\mbox{mod}\ 4). 
\end{array} \right.$$
Notice that in the first case we have $(12, s-6) = (s,6)$, while in the second case $(12,s-6) = 2(s,6)$.  Since $\Delta(-1) \cong 1_{\stmod \Lambda}$ if and only if $\chr(k)=2$, the stated periods for $\Lambda$ now follow. $\Box$\\

Finally, we conclude with a brief list of the remaining open cases.
\begin{itemize}
\item The periods of algebras of type $(\mathbb{D}_{n},s,2)$ with $n$ and $s$ even.  By \cite{Per} the period is either $2s(2n-3)/(s,n-1)$ or $4s(2n-3)/(s,n-1)$.
\item  The Calabi-Yau dimensions of algebras of type $(\mathbb{D}_n,s,2)$ with $n$ and $s$ even.  By 7.2, if $\stmod \Lambda$ is $d$-Calabi-Yau then $(s,n-1)=1$ and $d \equiv 1 - (n-1)^{-1}\ (\mbox{mod}\ 2s(2n-3))$.
\item The periods of algebras of type $(\mathbb{D}_4,s,3)$ when $3 | s$.  By \cite{Per} the period is either $5s$ or $15s$ in characteristic $2$ and either $10s/(s,2)$ or $30s/(s,2)$ otherwise.
\end{itemize}


\begin{thebibliography}{99}

\bibitem{Amiot} C. Amiot.  \emph{On the structure of triangulated categories with finitely many indecomposables.}   Bull. Soc. Math. France  135  (2007),  no. 3, 435--474.

\bibitem{DECSA} H. Asashiba.  \emph{The derived equivalence classification of representation-finite selfinjective algebras.}  J. Algebra 214 (1999), no. 1, 182--221.

\bibitem{Asa2} H. Asashiba.  \emph{On a lift of an individual stable equivalence to a standard derived equivalence for representation-finite self-injective algebras.}  Algebr. Represent. Theory 6 (2003), no. 4, 427--447.

\bibitem{Cover} H. Asashiba.  \emph{A generalization of Gabriel's Galois covering functors and derived equivalences.}  J. Algebra 334 (2011) 109--149.  


\bibitem{TEG} M. Auslander and I. Reiten.  \emph{On a theorem of E. Green on the dual of the transpose.}  Proc. ICRA V, CMS Conf. Proc. 11 (1991), 53-65.


\bibitem{DPA} J. Bia\l kowski, K. Erdmann and A. Skowro\'{n}ski.  \emph{Deformed preprojective algebras of generalized Dynkin type.}  Trans. Amer. Math. Soc. 359 (2007), no. 6, 2625-2650.

\bibitem{DMA} J. Bia\l kowski, K. Erdmann and A. Skowro\'{n}ski.  \emph{Deformed mesh algebras of generalized Dynkin type.} In preparation.

\bibitem{BiaSko} J. Bia\l kowski and A. Skowro\'{n}ski.  \emph{Calabi-Yau stable module categories of finite type.} Colloq. Math. 109 (2007), no. 2, 257-269.

\bibitem{GCY3} R. Bocklandt.  \emph{Graded Calabi-Yau algebras of dimension 3.} (Appendix by M. Van den Bergh) J. Pure Appl. Algebra 212 (2008), no. 1, 14--32.

\bibitem{PAAK} S. Brenner, M. C. R. Butler and A. D. King.  \emph{Periodic algebras which are almost Koszul.}  Algebras and Representation Theory  5 (2002), 331-367.

\bibitem{BMRRT} A. B. Buan, R. Marsh, M. Reineke, I. Reiten, and G. Todorov.  \emph{Tilting theory and cluster combinatorics.}  Adv. Math. 204 (2006), no. 2, 572--618.

\bibitem{Buch} R. O. Buchweitz.  \emph{Finite representation type and periodic Hochschild (co-)homology.}  Trends in the representation theory of finite-dimensional algebras (Seattle, WA, 1997),  81--109, Contemp. Math., 229, Amer. Math. Soc., Providence, RI, 1998. 


\bibitem{CiMa} C. Cibils and E. Marcos.  \emph{Skew category, Galois covering and smash product of a $k$-category.}  Proc. Amer. Math. Soc.  134  (2006),  no. 1, 39--50.


\bibitem{CoMo} M. Cohen and S. Montgomery.  \emph{Group-graded rings, smash products, and group actions.}  Trans. Amer. Math. Soc.  282  (1984),  no. 1, 237--258.

\bibitem{Die} E. Dieterich.  \emph{The Auslander-Reiten quiver of an isolated singularity.}  Singularities, representation of algebras, and vector bundles (Lambrecht, 1985),  244--264, Lecture Notes in Math., 1273, Springer, Berlin, 1987.

\bibitem{DW} E. Dieterich and A. Wiedemann.  \emph{The Auslander Reiten quiver of a simple curve singularity.} Trans. Amer. Math. Soc. 294 (1986), 455--475.

\bibitem{Per} A. Dugas.  \emph{Periodic resolutions and self-injective algebras of finite type.}  J. Pure Appl. Algebra 214 (2010), pp. 990-1000.  arXiv:0808.1311v2.

\bibitem{PHART} A. Dugas.  \emph{Periodicity of $d$-cluster-tilted algebras.}  Preprint (2010).  arXiv:1007.2811v1.



\bibitem{ErdSko}  K. Erdmann and A. Skowro\'nski.  \emph{The stable Calabi-Yau dimension of tame symmetric algebras.}  J. Math. Soc. Japan 58 (2006), 97-128.

\bibitem{ErdSko2} K. Erdmann and A. Skowro\'nski.  \emph{Periodic algebras.}  Trends in Representation Theory and Related Topics.  European Math. Soc. Series of Congress Reports, European Math. Soc. Publ. House, Zurich, 2008, 201--251.


\bibitem{ErdSna1} K. Erdmann and N. Snashall.  \emph{On Hochschild cohomology of preprojective algebras. I, II.} J. Algebra 205 (1998), no. 2, 391--412, 413--434.

\bibitem{CYFA} C.-H. Eu and T. Schedler.  \emph{Calabi-Yau Frobenius Algebras.}   J. Algebra  321  (2009),  no. 3, 774--815.

\bibitem{Green} E. L. Green.  \emph{Graphs with relations, coverings and group-graded algebras.}  Trans. Amer. Math. Soc.  279  (1983),  no. 1, 297--310.



\bibitem{DCFDA} D. Happel.  \emph{On the derived category of a finite-dimensional algebra.}  Comment. Math. Helv. 62 (1987), no. 3, 339-389.

\bibitem{TCRTA} D. Happel.  \emph{Triangulated categories in the representation theory of finite-dimensional algebras.} London Mathematical Society Lecture Note Series, 119. Cambridge University Press, Cambridge, 1988.

\bibitem{HCFDA} D. Happel.  \emph{Hochschild cohomology of finite-dimensional algebras.}  S\'eminaire d'Alg\`{e}bre Paul Dubreil et Marie-Paul Malliavin, 39\`{e}me Ann\'ee (Paris, 1987/1988),  108--126, Lecture Notes in Math., 1404, Springer, Berlin, 1989.


\bibitem{HPR1} D. Happel, U. Preiser and C. M. Ringel. \emph{Vinberg's characterization of Dynkin diagrams using subadditive functions with application to $D{\rm Tr}$-periodic modules.}  Representation theory, II (Proc. Second Internat. Conf., Carleton Univ., Ottawa, Ont., 1979), 280--294, Lecture Notes in Math., 832, Springer, Berlin, 1980.

\bibitem{HPR2} D. Happel, U. Preiser and C. M. Ringel.  \emph{Binary polyhedral groups and Euclidean diagrams.}  Manuscripta Math.  31  (1980), no. 1-3, 317--329.


\bibitem{HoJo1}  T. Holm and P. J\o rgensen.  \emph{Cluster categories, self-injective algebras and stable Calabi-Yau dimensions: type A.}  Preprint (2008) arXiv:math/0610728v2.

\bibitem{HoJo2}  T. Holm and P. J\o rgensen.  \emph{Cluster categories, self-injective algebras and stable Calabi-Yau dimensions: types D and E.}  Preprint (2008) arXiv:math/0612451v2.


\bibitem{TOC} B. Keller.  \emph{On triangulated orbit categories.}  Documenta Math. 10 (2005), 551-581.

\bibitem{TOCcor} B. Keller.  \emph{Corrections to `On triangulated orbit categories'.}  Available online at http://people.math.jussieu.fr/\~{ }keller/publ/index.html.

\bibitem{CYTC} B. Keller.  \emph{Calabi-Yau triangulated categories.}  Trends in Representation Theory of Algebras and Related Topics.  European Math. Soc. Series of Congress Reports, European Math. Soc. Publ. House, Zurich, 2008, 467--489.

\bibitem{CYTCcor} B. Keller.  \emph{Correction to `Calabi-Yau triangulated categories'.}  Available online at http://people.math.jussieu.fr/\~{ }keller/publ/index.html.

\bibitem{KelP} B. Keller.  \emph{On the Calabi-Yau property for higher cluster categories.}  Private communication (2009).

\bibitem{KeRe} B. Keller and I. Reiten.  \emph{Acyclic Calabi-Yau categories.} (Appendix by Michel Van den Bergh.)  Compos. Math.  144  (2008),  no. 5, 1332--1348.


\bibitem{DPic} J. Miyachi and A. Yekutieli.  \emph{Derived Picard groups of finite-dimensional hereditary algebras.}  Composito Mathematica 129 (2001), 341-368.

\bibitem{Noeth} I. Reiten and M. Van den Bergh.  \emph{Noetherian hereditary abelian categories satisfying Serre duality.}  J. Amer. Math. Soc. 15 (2002), no. 2, 295--366.

\bibitem{DCSE} J. Rickard.  \emph{Derived ategories and stable equivalence.}  J. Pure Appl. Algebra 61 (1989), no. 3, 303-317.

\bibitem{ADK} Ch. Riedtmann.  \emph{Algebren, Darstellungsk\"ocher,  Ueberlagerungen und zur\"uck.}  Comment. Math. Helvetici 55 (1980), 199--224.

\bibitem{Yosh} Y. Yoshino. \emph{Cohen-Macaulay modules over Cohen-Macaulay rings.}  London Mathematical Society Lecture Note Series, vol. 146,  Cambridge Univ. Press, Cambridge, 1990.

\end{thebibliography}
\end{document}